\documentstyle[]{article}

\title{The Rohlin Property for ${\bf Z}^{2}$-Actions on UHF Algebras}
\author{Hideki Nakamura}
\date{June\ \ 1996}

\begin{document}

\maketitle

\begin{abstract}
We define a Rohlin property for ${\bf Z}^{2}$-actions on 
UHF algebras and show a non-commutative Rohlin type theorem. 
Among those actions with the Rohlin property, 
we classify product type actions up to outer conjugacy. 
We consider two classes of UHF algebras. 
For UHF algebras in one class including the CAR algebra, 
there is one and only one outer conjugacy class of 
product type actions and for UHF algebras in the other class, 
contrary to the case of ${\bf Z}$-actions, 
there are infinitely many outer conjugacy classes of 
product type actions.
\end{abstract}

\vspace{4mm}

\newtheorem{definition}{Definition}
\newtheorem{lemma}[definition]{Lemma}
\newtheorem{proposition}[definition]{Proposition}
\newtheorem{theorem}[definition]{Theorem}
\newtheorem{remark}[definition]{Remark}
\newtheorem{corollary}[definition]{Corollary}

\newcommand{\qed}{\hbox{\rule[-2pt]{3pt}{6pt}}}

\section{Introduction}

A non-commtative Rohlin property was introduced by A. Connes for 
classification of (single) automorphisms
of von Neumann algebras (\cite{connes1,connes2}), 
and this property was generalized 
for example by A. Ocneanu 
(\cite{ocneanu1,ocneanu2}) to systems of commuting 
automorphisms and further to actions 
of discrete amenable groups.  
On the other hand, this notion has also proved useful in the framework 
of $C^{\ast}$-algebras 
(\cite{bek,eek,ek,ho1,ho2,kishimoto1,kishimoto2}), 
and in \cite{kishimoto1,kishimoto2} 
Kishimoto established a non-commutative Rohlin type 
theorem for automorphisms of 
UHF algebras (and some AF algebras) and classified automorphisms 
(i.e., ${\bf Z}$-actions) with the 
Rohlin property up to outer conjugacy.

The purpose of the paper is to extend Kishimoto's work to 
${\bf Z}^2$-actions.  Motivated by \cite{kishimoto1,kishimoto2},
in Section 2 we introduce notions of Rohlin property and uniform 
outerness to ${\bf Z}^N$-actions 
on unital $C^{\ast}$-algebras.  In the UHF algebra case and $N=1$, 
the uniform outerness was shown to be the same 
as the ordinary outerness of the relevant automorphism on 
the GNS von Neumann algebra obtained via the trace (\cite{kishimoto1}).  
Our main theorem here says that for ${\bf Z}^2$-actions 
on UHF algebras the Rohlin property 
characterizes the uniform outerness.  The main idea 
of the proof is similar to the one in \cite{kishimoto2}, 
but to avoid additional technical problems we make use 
of the stability i.e., the vanishing of $1$-cohomology 
obtained in \cite{ho1}.

In Section 3 we introduce three notions of conjugacy to 
${\bf Z}^2$-actions, i.e, approximate conjugacy, cocycle
conjugacy, and outer conjugacy.  Using the generalized 
determinant introduced by P. de la Harpe and G.
Skandalis (\cite{hs}), we show that approximate conjugacy 
implies cocycle conjugacy when a (unital) $C^{\ast}$-algebra
is simple and possesses a unique trace.
\par
In Section 4 we consider product type ${\bf Z}^2$-actions 
on UHF algebras, i.e., pairs $(\alpha, \beta)$ of
commuting automorphisms
$$
\alpha \cong \bigotimes_{k=1}^{\infty} Ad u_k, \ 
\beta  \cong \bigotimes_{k=1}^{\infty} Ad v_k.
$$
Considering the case when the $n_k \times n_k$ matrices 
$u_k, v_k$ commute at first, we show that the
Rohlin property in this case is characterized 
by the property of uniform distribution of the joint spectral 
set ${\rm Sp}(\otimes_{k=m}^n u_k, \otimes_{k=m}^n v_k) \ (m \leq n)$.  
From this we show that any such pairs are
approximately conjugate. 
\par
We then investigate two special classes of UHF algebras.  
The first one is of the form
$
\otimes_{k=1}^{\infty} M_{p_k^{i_k}},
$
where $p_k \ (k \in {\bf N})$ are primes and non-negative 
(finite) integers $i_k \ (k \in {\bf N})$ satisfy
$\sum_{k=1}^{\infty} i_k=\infty$.  The second one is of the form
$
\otimes_{k \in K} M_{q_k^{\infty}}
$
with primes $q_k \ (k \in K)$ (where $\# K \leq \infty$) 
and $M_{q_k^{\infty}}$ of course means the 
infinite tensor product of $M_{q_k}$.  For algebras 
in the first class we construct infinitely many
non-cocycle conjugate product type ${\bf Z}^2$-actions 
with the Rohlin property.  We would like to emphasize 
that for ${\bf Z}$-actions this phenomenon does not occur.  
On the other hand, for algebras in the
second class we show that all the product type 
${\bf Z}^2$-actions with the Rohlin property are
mutually approximately conjugate.  Combining these results 
we get the classification of product type
${\bf Z}^2$-actions with the Rohlin property on UHF algebras 
up to outer conjugacy.

\section{Rohlin type theorem}

Let $N$ be a positive integer . We first define the Rohlin property 
for ${\bf Z}^{N}$-actions on unital $C^{\ast}$-algebras. 
As mentioned above this is a simple generalization of that 
in the case of $N=1$ \cite{kishimoto1}. Let $\xi_{1},\ldots ,\xi_{N}$ 
be the canonical basis of ${\bf Z}^{N}$\ \ i.e.,
\[
\xi_{i}= (0,\ldots ,0,1,0,\ldots ,0),
\]
where $1$ is in the $i$-th component, and let $I= (1,\ldots ,1)$ 
throughout this section. For $m=(m_{1},\ldots ,m_{N})$ 
and $n=(n_{1},\ldots ,n_{N})\in {\bf Z}^{N}$,\ \ $m\leq n$ means 
$m_{i}\leq n_{i}$ for each $i=1,\ldots ,N$. We define 
\[
m{\bf Z}^{N}=
\{ (m_{1}n_{1},\ldots ,m_{N}n_{N})|
(n_{1},\ldots ,n_{N})\in{\bf Z}^{N} \}
\]
for $m=(m_{1},\ldots ,m_{N})\in{\bf Z}^{N}$ and 
let ${\bf Z}^{N}$ act on ${\bf Z}^{N}/m{\bf Z}^{N}$ by 
addition modulo $m{\bf Z}^{N}$.

\begin{definition}\label{lab-100}
{\rm 
Let $\alpha$ be a ${\bf Z}^{N}$-action 
on a unital $C^{\ast}$-algebra $A$\ \ 
i.e., $\alpha$ is a group homomorphism from ${\bf Z}^{N}$ 
into the automorphisms $Aut(A)$ of $A$. Then $\alpha$ 
is said to have the {\sl Rohlin property} if 
for any $m\in {\bf N}^{N}$ 
there exist $R\in {\bf N}$ and $m^{(1)},\ldots ,m^{(R)}\in {\bf N}^{N}$ 
with $m^{(1)},\ldots ,m^{(R)} \geq m$ 
and which satisfy the following condition: 
For any $\varepsilon >0$ and finite subset $F$ of $A$, 
there exist projections
\[
e^{(r)}_{g}\ \ 
\left(r=1,\ldots ,R,\ g\in {\bf Z}^{N}/m^{(r)}{\bf Z}^{N}\right)
\]
in $A$ satisfying
\[
\sum_{r=1}^{R}
\sum_{g\in {\bf Z}^{N}/m^{(r)}{\bf Z}^{N}}
e^{(r)}_{g}=1\ ,
\]
\begin{equation}
\| [x,\ e^{(r)}_{g}]\| <\varepsilon\ ,
\label{labcom}
\end{equation}
\[
\| \alpha_{{\xi}_{i}}(e^{(r)}_{g})-e^{(r)}_{{\xi}_{i}+g}\| <\varepsilon
\]
for any $x\in F$, $r=1,\ldots ,R$,  $i=1,\ldots ,N$,  
and $g\in {\bf Z}^{N}/m^{(r)}{\bf Z}^{N}$\ .
}
\end{definition}

\begin{remark}\label{lab-99}
{\rm 
When $A$ is a UHF algebra, 
using Christensen's perturbation argument 
(\cite[Theorem\ 5.3.]{christensen}),  we can restate 
the definition of the Rohlin property as follows.  
For any $n,m\in {\bf N}$ with $1\leq n\leq N$ 
there exist $R\in {\bf N}$ and 
positive integers $m^{(1)},\ldots ,m^{(R)} \geq m$ 
which satisfy the following condition: 
For any $\varepsilon >0$ and finite subset $F$ of $A$ 
there exist projections
\[
e^{(r)}_{0},\ldots ,e^{(r)}_{m^{(r)}-1}\ \ (r=1,\ldots ,R)
\]
in $A$ satisfying
\[
\sum_{r=1}^{R}
\sum_{j=0}^{m^{(r)}-1}
e^{(r)}_{j}=1\ ,
\]
\[
\| [x,\ e^{(r)}_{j}]\| <\varepsilon
\]
for each $r=1,\ldots ,r,\ j=0,\ldots ,m^{(r)}-1$ and $x\in F$, and
\[
\| \alpha_{\xi_{n}}(e^{(r)}_{j})-e^{(r)}_{j+1}\| <\varepsilon\ ,
\]
\[
\| \alpha_{\xi_{n'}}(e^{(r)}_{j})-e^{(r)}_{j}\| <\varepsilon
\]
for each $n'=1,\ldots ,N$ with $n'\neq n,\ \ r=1,\ldots R$ 
and $j=0,\ldots ,m^{(r)}-1$, 
where $e^{(r)}_{m^{(r)}}\equiv e^{(r)}_{0}$ .
}
\end{remark}

For automorphisms of $C^{\ast}$-algebras a notion of 
uniform outerness was introduced in \cite{kishimoto1}. 
That is, an automorphism $\alpha$ of a unital 
$C^{\ast}$-algebra $A$ 
is said to be {\sl uniformly outer} if 
for any $a\in A$, any nonzero projection $p\in A$ and 
any $\varepsilon >0$,  
there exist projections $p_{1},\ldots ,p_{n}$ in $A$ such that 
\[
p=\sum_{i=1}^{n}p_{i}\ ,
\]
\[
\| p_{i}a\alpha(p_{i})\| <\varepsilon\ \ (i=1,\ldots ,n)\ .
\]
It was shown that 
this notion for automorphisms of UHF algebras 
is equivalent to the usual outerness 
for the automorphisms of the von Neumann algebras 
obtained through the GNS representations associated 
with the traces (\cite[Theorem\ 4.5]{kishimoto1}).
Based on this fact and the Rohlin type theorem for 
automorphisms of von Neumann algebras due to 
A. Connes (\cite[Theorem\ 1.2.5]{connes2}), 
a $C^{\ast}$-algebraic version of the theorem 
(for the UHF algebras) was shown by A. Kishimoto 
(\cite[Theorem\ 1.3]{kishimoto2}). We extend Kishimoto's 
work to ${\bf Z}^{2}$-actions.

\begin{theorem}\label{lab1}
Let $\alpha$ be a ${\bf Z}^{2}$-action on a UHF algebra $A$. 
Then the following conditions are equivalent:
\begin{list}{}{}
\item[{\rm (1)}] $\alpha$ has the Rohlin property.
\item[{\rm (2)}] ${\alpha}_{g}$ is uniformly outer 
for each $g\in {\bf Z}^{2}\setminus \{0\}$.
\end{list}
\end{theorem}

Once we establish this theorem, we have immediately

\begin{corollary}\label{lab23}
Let $\alpha$ be a ${\bf Z}^{2}$-action on a UHF algebra $A$. 
Then the following conditions are equivalent:
\begin{list}{}{}
\item[{\rm (1)}] $\alpha$ has the Rohlin property as a 
${\bf Z}^{2}$-action 
on $A$.
\item[{\rm (2)}] ${\alpha}_{g}$ has the Rohlin property 
as an automorphism of $A$ for each $g\in {\bf Z}^{2}\setminus \{0\}$.
\end{list}
\end{corollary}

In Theorem \ref{lab1} it is obvious that (1) implies (2). 
We devote the rest of this section to prove the converse 
in several steps.

\begin{lemma}\label{lab2}
Let $\alpha$ be a ${\bf Z}^{2}$-action on a UHF algebra $A$. 
If the condition {\rm (2)} in Theorem \ref{lab1} holds 
then for any $m=(m_{1},m_{2})\in {\bf N}^{2}$, $\varepsilon >0$ 
and any unital full matrix subalgebra $B$ of $A$, 
there exists an orthogonal family 
$(\ e_{g}\ |\ g\in {\bf Z}^{2},0\leq g\leq m-I\ )$ 
of projections in $A\cap B'$ such that
\[
\|{\alpha}_{{\xi}_{i}}(e_{g})-e_{{\xi}_{i}+g}\|<\varepsilon
\]
for any $i=1,\ 2$ and $g\in {\bf Z}^{2}$ 
with $\ 0\leq g,{\xi}_{i}+g\leq m-I$, and furthermore,
\[
1<(|m|+1)\tau (e_{0})\ ,
\]
where $\tau$ is the unique tracial state of $A$ 
and $|m|\equiv m_{1}\cdot m_{2}$.
\end{lemma}

\noindent $Proof.$\ \ Let $({\pi}_{\tau},H_{\tau})$ be 
the GNS representation associated with $\tau$. 
By the uniqueness of a trace we can extend each 
${\alpha}_{g}\ (g\in {\bf Z}^{2})$ to an automorphism 
of the AFD ${\rm II_{1}}$ factor 
${\pi}_{\tau}(A)''\ (\subseteq B(H_{\tau}))$ 
and we use the same symbol ${\alpha}_{g}$ 
for this extension. Since ${\alpha}_{g}$ is 
outer on ${\pi}_{\tau}(A)''$ for  
$g\in {\bf Z}^{2}\setminus \{0\}$ by 
\cite[Theorem\ 4.5]{kishimoto1}, 
it follows from \cite[Theorem\ 2]{ocneanu1} that for any 
$m\in {\bf N}^{2}$ there exists a strongly central sequence 
\[
\left( (\, E^{(j)}_{g}\, |\, g\in {\bf Z}^{2},0\leq g\leq m-I\, )\ 
| \ j\in{\bf N} \right)\
\]
of orthogonal families of projections in ${\pi}_{\tau}(A)''$ 
such that
\[
\sum_{\stackrel{\scriptstyle g\in {\bf Z}^{2}}{0\leq g\leq m-I}}
E^{(j)}_{g}=1
\]
for each $j\in {\bf N}$ and
\[
{\alpha}_{{\xi}_{i}}(E^{(j)}_{g})-E^{(j)}_{{\xi}_{i}+g}
\longrightarrow 0
\]
strongly as $j\rightarrow\infty$ for each $i=1,\ 2$ 
and $g\in {\bf Z}^{2}$ with  $0\leq g\leq m-I$, 
where $(m_{1},k)\equiv (0,k),\ (k,m_{2})\equiv (k,0)$.

From this central sequence we shall construct 
a uniformly central sequence  
\[
\left( (\, f^{(j)}_{g}\, |\, g\in {\bf Z}^{2},\ 0\leq g\leq m-I\, )\ 
|\ j\in{\bf N}\right)
\]
of orthogonal families of projections in $A$ such that
\[
{\pi}_{\tau}\left(\sum_{\stackrel{\scriptstyle g\in {\bf Z}^{2}}
{0\leq g\leq m-I}}f^{(j)}_{g}\right)\longrightarrow 1
\]
strongly as $j\rightarrow\infty$ and
\[
\| {\alpha}_{{\xi}_{i}}(f^{(j)}_{g})-f^{(j)}_{{\xi}_{i}+g}\|
\longrightarrow 0
\]
as $j\rightarrow\infty$ for each $i=1,\ 2$ and $g\in {\bf Z}^{2}$ 
with $0\leq g,g+\xi_{i}\leq m-I$. To do this, 
let $(\, A_{j}\, |\, j\in {\bf N}\, )$ be an increasing sequence 
of unital full matrix subalgebras of $A$ such that $\cup_{j}A_{j}$ 
is dense in $A$. From \cite[Lemma\ 4.7]{kishimoto1} 
we find a uniformly central sequence 
$(\, e_{j}\,|\, j\in {\bf N}\, )$ of projections in $A$ 
such that
\[
{\pi}_{\tau}(e_{j})-E^{(j)}_{0}\longrightarrow 0
\]
strongly as $j\rightarrow\infty$. Changing $e_{j}$ slightly 
and taking a subsequence, we may assume 
that $e_{j}\in (\cup_{k}A_{k})\cap A_{j}'$ 
for each $j$. Let $\varepsilon >0$. 
From \cite[Corollary\ 6.8]{christensen}, 
by taking inner perturbations, 
there are $\alpha_{1},\alpha_{2}\in Aut(A)$ such that
\[
\| \alpha_{i}-\alpha_{\xi_{i}}\| <\varepsilon\ ,
\]
\[
\alpha_{i}\left(\bigcup_{j\in{\bf N}}A_{j}\right)
\subseteq \bigcup_{j\in{\bf N}}A_{j}
\]
for $i=1,2$. Set
\[
h_{j}=e_{j}\left(
\sum_{\stackrel{\scriptstyle g=(g_{1},g_{2})
\in {\bf Z}^{2}\setminus\{ 0\}}{-(m-I)\leq g\leq m-I}}
({\alpha}_{1})^{g_{1}}
({\alpha}_{2})^{g_{2}}(e_{j})\right) e_{j}\ ,
\]
\[
k_{j}=e_{j}\left(
\sum_{\stackrel{\scriptstyle g\in {\bf Z}^{2}\setminus\{ 0\}}
{-(m-I)\leq g\leq m-I}}{\alpha}_{g}(e_{j})\right) e_{j}\ ,
\]
$\eta_{j}=\tau (h_{j})$ and $\kappa_{j}=\tau (k_{j})$. 
Then it follows that
\[
\| h_{j}-k_{j}\| <\varepsilon (m_{1}+m_{2}-2)
\{ (2m_{1}-1)(2m_{2}-1)-1\}\ .
\]
Furthermore, $\lim_{j}\kappa_{j}=0$ since
\begin{eqnarray*}
\lim_{j\rightarrow\infty}{\kappa}_{j}
  &=& \lim_{j\rightarrow\infty}\tau
\left( e_{j}\sum_{\stackrel{\scriptstyle g
\in {\bf Z}^{2}\setminus\{ 0\}}{-(m-I)\leq g\leq m-I}}
{\alpha}_{g}(e_{j})\right)\\
  & \leq & \lim_{j\rightarrow\infty}
\sum_{\stackrel{\scriptstyle g,\ h\in {\bf Z}^{2},\ g\neq h}
{0\leq g,\ h\leq m-I}}\tau
(e_{j}{\alpha}_{g-h}(e_{j}))\\
  &\leq & \lim_{j\rightarrow\infty}
\sum_{\stackrel{\scriptstyle g,\ h\in {\bf Z}^{2},\ g\neq h}
{0\leq g,\ h\leq m-I}}\tau
({\alpha}_{h}(e_{j}){\alpha}_{g}(e_{j}))\\
  &=&\lim_{j\rightarrow\infty}
\sum_{\stackrel{\scriptstyle g,\ h\in {\bf Z}^{2},\ g\neq h}
{0\leq g,\ h\leq m-I}}\tau (E^{(j)}_{h}E^{(j)}_{g})=0\ .
\end{eqnarray*}
Let $p_{j}$ be the spectral projection 
of $h_{j}$ corresponding to 
$(0,\ {\eta}_{j}^{\frac{1}{2}})$. 
Then $p_{j}\in A$ since ${\rm Sp}(h_{j})$ 
is finite. $p_{j}\leq e_{j}$ and 
${{\eta}_{j}}^{\frac{1}{2}}(e_{j}-p_{j})\leq h_{j}$ 
because ${\eta_{j}}^{\frac{1}{2}}
{\chi}_{[{{\eta}_{j}}^{\frac{1}{2}},\ \infty )}
(t)\leq t$ ($t\in [0,\ \infty))$, hence
\[
\tau (e_{j})-{{\eta}_{j}}^{\frac{1}{2}}\leq 
\tau (p_{j})\leq \tau (e_{j})\ .
\]
In addition
\[
\| p_{j}\left( 
\sum_{\stackrel{\scriptstyle g=(g_{1},g_{2})
\in {\bf Z}^{2}\setminus\{ 0\}}
{-(m-I)\leq g\leq m-I}}({\alpha}_{1})^{g_{1}}
({\alpha}_{2})^{g_{2}}(p_{j})\right)p_{j}\| \leq
\| p_{j}h_{j}p_{j}\| \leq {{\eta}_{j}}^{\frac{1}{2}}.
\]
So for any $g,h\in {\bf Z}^{2}$ 
with $0\leq g,h\leq m-I$ and $g\neq h$, we have
\begin{eqnarray*}
\| {\alpha}_{g}(p_{j}){\alpha}_{h}(p_{j}){\|}^{2}
  &=&\| {\alpha}_{g}(p_{j}{\alpha}_{h-g}(p_{j})){\|}^{2}
=\| p_{j}{\alpha}_{h-g}(p_{j}){\|}^{2}\\
  &=&\| p_{j}{\alpha}_{h-g}(p_{j})p_{j}\| \\
  &\leq&\| p_{j}
\left( \sum_{\stackrel{\scriptstyle g\in {\bf Z}^{2}\setminus\{ 0\}}
{-(m-I)\leq g\leq m-I}}{\alpha}_{g}(p_{j})\right)p_{j}\| \\
  &\leq&\varepsilon '+{{\eta}_{j}}^{\frac{1}{2}},
\end{eqnarray*}
where $\varepsilon '= \varepsilon 
(m_{1}+m_{2}-2)\{(2m_{1}-1)(2m_{2}-1)-1\}$. 
Here $\lim\eta_{j}=0$ since 
$\lim\kappa_{j}=0$ and 
$\lim\| h_{j}-k_{j}\|=0$.
Therefore taking a sufficiently large $j$ 
for each $\varepsilon >0$, 
we obtain the required $f^{(j)}_{g}$ 
near ${\alpha}_{g}(p_{j})$ by slight modification.

Noting that 
$\sum_{0\leq g\leq m-I}\tau(f^{(j)}_{g})\rightarrow 1$ 
and $\tau (f^{(j)}_{g})=\tau (f^{(j)}_{h})$, we have
\[
\tau (f^{(j)}_{g})\longrightarrow \frac{1}{|m|}\ .
\]
Furthermore for any unital full matrix subalgebra 
$B$ of $A$, taking a sufficiently large $j$, 
we may assume that $f^{(j)}_{g}\in A\cap B'$ 
for any $g$. This concludes the proof. \hfill \qed

\vspace{4mm}

In Ocneanu's result \cite[Theorem\ 2]{ocneanu1}, 
applied in the above proof, we have the cyclicity condition 
(under the action) of the projections 
$(\, E^{(j)}_{g}\, |\, g\in {\bf Z}^{2},0\leq g\leq m-I\, )$ 
in the von Neumann algebra $\pi_{\tau}(A)''$. 
However when approximating these projections by the projections 
$(\, f^{(j)}_{g}\, |\, g\in {\bf Z}^{2},0\leq g\leq m-I\, )$ 
in the $C^{\ast}$-algebra $A$, 
we lose the cyclicity condition. 
It is our next problem to restore this cyclicity condition. 
To do this we need a technical lemma 
from \cite{kishimoto2}. 
Let $K(l^{2}({\bf Z}))$ 
be the compact operators on $l^{2}({\bf Z})$ 
and let $(\, E_{i,j}\, |\, i,j\in{\bf Z}\, )$ 
be the canonical matrix units for $K(l^{2}({\bf Z}))$. 
On $K(l^{2}({\bf Z}))$ we define an automorphism $\sigma$ 
by $\sigma (E_{i,j})=E_{i+1,j+1}\ (i,\ j\in {\bf Z})$. 
For any $n,k,l\in {\bf N}$ with $1<k<l$, define 
\[
N=n(2k+l-1)\ ,
\]
\[
f=\sum_{i=1}^{k-1} 
\left( \frac{i}{k}E_{ni,ni}+\frac{k-i}{k}E_{n(k+l+i),n(k+l+i)}
+\frac{\sqrt{i(k-i)}}{k}E_{ni,n(k+l+i)}\right.
\]
\begin{equation}
+\left.\frac{\sqrt{i(k-i)}}{k}E_{n(k+l+i),ni}\right)
+\sum_{i=k}^{k+l}E_{ni,ni}\ ,
\label{hoshi}
\end{equation}
\[
e_{i}=\sigma^{i-n}(f)\ \ (i=0,\ldots ,n-1)\ .
\]
Then
$(\, e_{i}\, |\, i=0,\ldots ,n-1\, )$ 
is an orthogonal family of projections in $K(l^{2}({\bf Z}))$. 
Hence for any $\varepsilon >0$, 
there exist $k,l$ with $1\ll k\ll l$ such that 
\[
\sum_{i=0}^{n-1}e_{i}\leq P_{N}\ \ 
(P_{N}\equiv\sum_{i=0}^{N-1}E_{i,i})\ ,
\]
\[
\| \sigma (e_{i})-e_{i+1}\| <\varepsilon\ \ 
(i=0,\ldots ,n-1,\ \ e_{n}\equiv e_{0})\ ,
\]
\[
\frac{n\ {\rm dim}\ e_{0}}{N}>1-\varepsilon 
\]
(see \cite[Lemma\ 2.1]{kishimoto2} for the detail).
Using these estimates we have the next lemma.

\begin{lemma}\label{lab4}
Let $\alpha$ be a ${\bf Z}^{2}$-action on a UHF algebra $A$. 
If $\alpha_{g}$ is uniformly outer for any 
$g\in{\bf Z}^{2}\setminus \{ 0\}$, 
then for any $m\in {\bf N}$, $\varepsilon >0$ 
and any unital full matrix subalgebra $B$ of $A$ 
there exists an orthogonal family 
$(\, e_{i}\, |\, i=0,\ldots ,m-1\, )$ 
of projections in $A\cap B'$ such that
\[
\|\alpha_{\xi_{1}}(e_{i})-e_{i+1}\| <\varepsilon\ ,
\]
\[
\|\alpha_{\xi_{2}}(e_{i})-e_{i}\| <\varepsilon\ ,
\]
\[
\tau (1-\sum_{i=0}^{m-1}e_{i})\leq \varepsilon\tau (e_{0})
\]
for $i=0,\ldots ,m-1$, where $e_{m}\equiv e_{0}$.
\end{lemma}

\noindent {\it Proof}.\ \ Let $m\in {\bf N}$, 
$\varepsilon_{1}>0$ and let $B_{1}$ be 
a unital full matrix subalgebra of $A$. 
By the above statement there exist 
$k_{1},\ l_{1}\in {\bf N}$ with $1\ll k_{1}\ll l_{1}$ 
and an orthogonal family 
$(\, e_{i}\, |\, i=0,\ldots ,m-1\, )$ 
of projections in $K(l^{2}({\bf Z}))$ such that
\[
\sum_{i=0}^{m-1}e_{i}\leq P_{N_{1}}\ ,
\]
\[
\| \sigma (e_{i})-e_{i+1} \|
<\varepsilon_{1}\ \ (i=0,\ldots ,m-1)\ ,
\]
\begin{equation}
\frac{m\ {\rm dim}\ e_{0}}{N_{1}}
>1-\varepsilon_{1}\ ,
\label{1star}
\end{equation}
where $N_{1}\equiv m(2k_{1}+l_{1}-1)$  
and $e_{m}\equiv e_{0}$. 
Similarly by the above statement   
(for $n=1,\ \varepsilon_{1}$ and $B_{1}$), 
there exist $k_{2},\ l_{2}\in {\bf N}$ 
with $1\ll k_{2}\ll l_{2}$ 
and a projection $e$ in $K(l^{2}({\bf Z}))$ such that
\[
e\leq P_{N_{2}}\ ,
\]
\[
\| \sigma (e)-e \|<\varepsilon_{1}\ ,
\]
\begin{equation}
\frac{{\rm dim}\ e}{N_{2}}>1-\varepsilon_{1}\ ,
\label{laby}
\end{equation}
where $N_{2}\equiv 2k_{2}+l_{2}-1$.

Next by applying Lemma \ref{lab2} to  
$(N_{1},N_{2})\in {\bf N}^{2}$, 
any $\varepsilon_{2}>0$ and $B_{1}$, 
there exists an orthogonal family 
$(\, p_{g}\, |\, g\in {\bf Z}^{2},0\leq g\leq (N_{1}-1,N_{2}-1)\, )$ 
of projections in $A\cap B_{1}'$ such that
\[
\| \alpha_{\xi_{i}}(p_{g})-p_{\xi_{i}+g}\| 
<\varepsilon_{2}\ ,
\]
\begin{equation}
1\leq (N_{1}N_{2}+1)\tau (p_{0})
\label{labz}
\end{equation}
for any $i=1,\ 2$ and $g\in {\bf Z}^{2}$ 
with $0\leq g,\xi_{i}+g\leq (N_{1}-1,N_{2}-1)$.

If we put
\[
x_{1}=\frac{1}{N_{2}}\sum_{j=0}^{N_{2}-1} 
\left\{\sum_{i=0}^{N_{1}-2}p_{(i+1,j)}
\alpha_{\xi_{1}}(p_{(i,j)})\right.
\]
\[
+\left.(1-\sum_{i=1}^{N_{1}-1}p_{(i,j)})
(1-\sum_{i=0}^{N_{1}-2}\alpha_{\xi_{1}}(p_{(i,j)}))\right\}\ , 
\]
then we have 
\[
x_{1}\alpha_{\xi_{1}}(p_{(i,j)})=p_{(i+1,j)}x_{1}
\]
for any $i=0,\ldots ,N_{1}-2,\ j=0,\ldots ,N_{2}-1$ 
and
\[
x_{1}-1= \frac{1}{N_{2}}\sum_{j=0}^{N_{2}-1} 
\left\{\sum_{i=0}^{N_{1}-2}p_{(i+1,j)}
\left( \alpha_{\xi_{1}}(p_{(i,j)})-p_{(i+1,j)}\right)\right.
\]
\[
+\left. (1-\sum_{i=1}^{N_{1}-1}p_{(i,j)})
(-\sum_{i=0}^{N_{1}-2}\alpha_{\xi_{1}}(p_{(i,j)}))\right\}\ .
\]
Noting that 
$\| \alpha_{\xi_{1}}(p_{(i,j)})-p_{(i+1,j)}\| 
<\varepsilon_{2}$, 
we have $\| x_{1}-1\| <2(N_{1}-1)\varepsilon_{2}$. 
So taking the polar decomposition $u_{1}|x_{1}|$ of $x_{1}$ 
for a sufficiently small $\varepsilon_{2}>0$, 
we obtain a unitary $u_{1}$ with $\| u_{1}-1\| 
<4(N_{1}-1)\varepsilon_{2}$. By the uniqueness 
of the polar decomposition we have
\[
Ad\, u_{1}\circ\alpha_{\xi_{1}}(p_{(i,j)})=p_{(i+1,j)}
\]
for $i=0,\ldots ,N_{1}-2,\ j=0,\ldots ,N_{2}-1$. 
Similarly for $\alpha_{\xi_{2}}$, 
we obtain a unitary $u_{2}$ in $A$ such that 
\[
\| u_{2}-1\| < 4(N_{2}-1)\varepsilon_{2}\ ,
\]
\[Ad\, u_{2}\circ\alpha_{\xi_{2}}(p_{(i,j)})=p_{(i,j+1)}
\]
for $i=0,\ldots ,N_{1}-1,\ j=0,\ldots ,N_{2}-2$. 
Let $\alpha_{1}= Ad\, u_{1}\circ\alpha_{\xi_{1}}$ 
and let $\alpha_{2}= Ad\, u_{2}\circ\alpha_{\xi_{2}}$. 
Since $[p_{(0,0)}]=[p_{(1,0)}]$ 
it follows that there exists a partial isometry 
$v_{1}$ of $A\cap B_{1}'$ such that 
$v_{1}^{\ast}v_{1}=p_{(0,0)}$ 
and $v_{1}v_{1}^{\ast}=p_{(1,0)}$. 
Similarily there exists a partial isometry $v_{2}$ 
of $A\cap B_{1}'$ such that $v_{2}^{\ast}v_{2}=p_{(0,0)}$ 
and $v_{2}v_{2}^{\ast}=p_{(0,1)}$. 
Then $Ad\, v_{2}^{\ast}\circ\alpha_{2}(p_{(0,0)})=p_{(0,0)}$, 
so $Ad\, v_{2}^{\ast}\circ\alpha_{2}\in 
Aut(p_{(0,0)}Ap_{(0,0)})$. 
On the other hand $\alpha_{\xi_{2}}\in Aut(A)$ 
has the Rohlin property as a single automorphism, 
and hence so does $Ad\, v_{2}^{\ast}\circ\alpha_{2}$ . 
Therefore 
$Ad\, v_{2}^{\ast}\circ\alpha_{2}$ is  
stable by \cite{ho1,ek}. 
More precisely for any $\varepsilon_{3}>0$, 
any unital full matrix subalgebra $B_{2}$ 
of $A$ and the unitary 
$v_{2}^{\ast}\alpha_{2}(v_{1})^{\ast}\alpha_{1}(v_{2})v_{1}
\in p_{(0,0)}Ap_{(0,0)}$, 
if $B_{1}$ is taken sufficiently large  
in advance, we have a unitary $w$ in $A\cap B_{2}'$ such that
\[
\| v_{2}^{\ast}\alpha_{2}(v_{1})^{\ast}
\alpha_{1}(v_{2})v_{1}-(Ad\, v_{2}^{\ast}
\circ\alpha_{2}(w))\cdot w^{\ast}\| <\varepsilon_{3}.
\]
Let $w_{1}= v_{1}w$ and let $w_{2}= v_{2}$. 
Then $w_{1}$ and $w_{2}$ are partial isometries 
in $A\cap B_{2}'$ such that $w_{1}^{\ast}w_{1}
=w_{2}^{\ast}w_{2}=p_{(0,0)}$, $w_{1}w_{1}^{\ast}
=p_{(1,0)}$, $w_{2}w_{2}^{\ast}=p_{(0,1)}$ and
\begin{equation}
\| \alpha_{1}(w_{2})w_{1}-\alpha_{2}(w_{1})w_{2}\| =
\| v_{2}^{\ast}\alpha_{2}(v_{1})^{\ast}(\alpha_{1}
(v_{2})v_{1}w-\alpha_{2}(v_{1}w)v_{2})w^{\ast}\|
\leq \varepsilon_{3}
\label{2star}
\end{equation}

Define
\[
E^{(k)}_{i,j}=\left\{
\begin{array}{ll}
\alpha_{2}^{i-1}(\alpha_{1}^{k}(w_{2})) 
\alpha_{2}^{i-2}(\alpha_{1}^{k}(w_{2})) 
\cdots \alpha_{2}^{j}(\alpha_{1}^{k}(w_{2})) 
& (\, i>j\, )\\
p_{(k,i)} & (\, i=j\, )\\
\alpha_{2}^{i}(\alpha_{1}^{k}(w_{2}))^{\ast} 
\alpha_{2}^{i+1}(\alpha_{1}^{k}(w_{2}))^{\ast} 
\cdots \alpha_{2}^{j-1}(\alpha_{1}^{k}(w_{2}))^{\ast} 
& (\, i<j\, )
\end{array}
\right. 
\]
for $k=0,1,\ i,j=0,\ldots ,N_{2}-1$. 
Then we can easily see that 
$(\, E^{(k)}_{i,j}\, |\, i,j=0,\ldots ,N_{2}-1\, )$ 
is a system of matrix units. 
For any unital full matrix subalgebra $B_{3}$ of $A$, 
by taking a sufficiently large $B_{2}$ including $B_{3}$, 
we may assume that 
$\{\, E^{(k)}_{i,j}\, |\, k=0,1;\ i,j=0,\ldots ,N_{2}-1\, \}
\subseteq A\cap B_{3}'$. Let $C^{(k)}$ be 
the $C^{\ast}$-subalgebra of $A$ generated by 
$\{\, E^{(k)}_{i,j}\, |\, i,j=0,\ldots ,N_{2}-1\, \}$ 
and let $\Phi_{k}$ be the canonical isomorphism 
from $C^{(k)}$ onto $P_{N_{2}}K(l^{2}({\bf Z}))P_{N_{2}}$. 
Define
\[
e^{(k)}=\Phi_{k}^{-1}(e)\ .
\]
Since 
$\sigma\circ\Phi_{k}=\Phi_{k}\circ\alpha_{2}\lceil C^{(k)}$, 
we have from (\ref{laby}) that 
\[
e^{(k)}\leq\sum_{i=0}^{N_{2}-1}p_{(k,i)}\ ,
\]
\[
\| \alpha_{2}(e^{(k)})-e^{(k)}\| <\varepsilon_{1}\ ,
\]
\begin{equation}
\tau (e^{(k)})>(1-\varepsilon_{1})N_{2}\tau (p_{(0,0)})\ .
\label{labzz}
\end{equation}
Furthermore define
\[
W_{1}=\left( \sum_{i=0}^{N_{2}-1}\alpha_{2}^{i}
(w_{1})\right) e^{(0)}\ ,
\]
\[
e^{(1)\prime}=W_{1}W_{1}^{\ast}
\leq\sum_{i=0}^{N_{2}-1}p_{(1,i)}\ .
\]
Again for any unital full matrix subalgebra 
$B_{4}$ of $A$, by taking a sufficiently large $B_{3}$ 
including $B_{4}$, we may assume that 
$W_{1}\in A\cap B_{4}'$. 
Then recalling the formula (\ref{hoshi}), 
we have
\begin{equation}
e^{(1)\prime}=
\left( \sum_{i=0}^{N_{2}-1}\alpha_{2}^{i}(w_{1})\right)
X
\left( \sum_{i'=0}^{N_{2}-1}\alpha_{2}^{i'}(w_{1})\right)^{\ast}\ ,
\label{hoshi2}
\end{equation}
where $X$ denotes
\[
\sum_{j=1}^{k_{2}-1} \left\{ 
\frac{j}{k_{2}}E^{(0)}_{j,j}
+\frac{k_{2}-j}{k_{2}}E^{(0)}_{k_{2}+l_{2}+j,k_{2}+l_{2}+j}
+\frac{\sqrt{j(k_{2}-j)}}{k_{2}}E^{(0)}_{j,k_{2}+l_{2}+j}
\right.
\]
\[
\left. 
+\frac{\sqrt{j(k_{2}-j)}}{k_{2}}E^{(0)}_{k_{2}+l_{2}+j,j}\right\}
+\sum_{j=k_{2}}^{k_{2}+l_{2}}E^{(0)}_{j,j}\ .
\]
On the right hand side of (\ref{hoshi2}), 
the nonzero terms are calculated as follows:
\begin{eqnarray*}
\alpha_{2}^{j}(w_{1})E^{(0)}_{j,j}
\alpha_{2}^{j}(w_{1})^{\ast}
&=&\alpha_{2}^{j}(w_{1})p_{(0,j)}
\alpha_{2}^{j}(w_{1})^{\ast}\\
&=&p_{(1,j)}\\
&=&E^{(1)}_{j,j}\ ,
\end{eqnarray*}
\begin{eqnarray*}
\lefteqn{\alpha_{2}^{j}(w_{1})E^{(0)}_{j,k_{2}+l_{2}+j}
\alpha_{2}^{k_{2}+l_{2}+j}(w_{1})^{\ast}}\\
&=&\alpha_{2}^{j}(w_{1})
\alpha_{2}^{j}(w_{2})^{\ast}\cdots
\alpha_{2}^{k_{2}+l_{2}+j-1}(w_{2}^{\ast}\alpha_{2}(w_{1})^{\ast})\\
&\stackrel{\varepsilon_{3}}{\approx}&
\alpha_{2}^{j}(w_{1})
\alpha_{2}^{j}(w_{2})^{\ast}\cdots
\alpha_{2}^{k_{2}+l_{2}+j-1}(w_{1}^{\ast}\alpha_{1}(w_{2})^{\ast})\ ,
\end{eqnarray*}
where $x\stackrel{\varepsilon}{\approx} y$ 
means $\| x-y\| <\varepsilon$. Applying (\ref{2star}) 
repeatedly we have
\[
\alpha_{2}^{j}(w_{1})E^{(0)}_{j,k_{2}+l_{2}+j}
\alpha_{2}^{k_{2}+l_{2}+j}(w_{1})^{\ast}
\]
\[
\stackrel{(k_{2}+l_{2})\varepsilon_{3}}{\approx}
E^{(1)}_{j,k_{2}+l_{2}+j}\ .
\]
We estimate the other nonzero terms similarly and obtain
\[
\| e^{(1)\prime} -e^{(1)}\| \leq
\sum_{j=1}^{k_{2}-1}\frac{\sqrt{j(k_{2}-j)}}{k_{2}}
(k_{2}+l_{2})\varepsilon_{3}\ .
\]
Let $c_{1}(k_{2},l_{2},\varepsilon_{3})$ 
be the right hand side of the above inequality. 
Then we have
\begin{eqnarray*}
\| \alpha_{2}(e^{(1)\prime} )-e^{(1)}\|
&\leq&2\| e^{(1)\prime} -e^{(1)}\|+\| 
\alpha_{2}(e^{(1)})-e^{(1)}\|\\
&\leq&2c_{1}(k_{2},l_{2},\varepsilon_{3})+\varepsilon_{1}\ .
\end{eqnarray*}
For any $\varepsilon_{4}>0$ 
and any unital full matrix subalgebra $B_{5}$ of $A$, 
applying \cite[Lemma\ 3.5]{kishimoto1} 
and taking a sufficiently large $B_{4}$ 
including $B_{5}$, we have a partial isometry 
$W_{1}'$ of $A\cap B_{5}'$ such that $(W_{1}' )^{\ast}W_{1}' 
=e^{(0)}$, $W_{1}' (W_{1}')^{\ast}=e^{(1)\prime }$ 
and
\begin{eqnarray}
\| \alpha_{2}(W_{1}')-W_{1}'\| 
&\leq&\| \alpha_{2}(e^{(0)})-e^{(0)}\| +
\| \alpha_{2}(e^{(1)\prime} )-e^{(1)\prime}\| 
+\varepsilon_{4}\nonumber \\
&\leq&\varepsilon_{1}+2c_{1}(k_{2},l_{2},\varepsilon_{3})
+\varepsilon_{1}+
\varepsilon_{4}\ .
\label{labw'}
\end{eqnarray}
Of course we can make the last quantity very small. 
By using this $W_{1}'$, let $D$ be the $C^{\ast}$-subalgebra 
of $A$ generated by 
\[
\{ \alpha_{1}^{i-1}(W_{1} )
\alpha_{1}^{i-2}(W_{1}')\cdots\alpha_{1}^{j}(W_{1}' )\ 
|\ N_{1}-1\geq i>j\geq 0\ \}\ .
\] 
Here again for any unital full matrix subalgebra $B_{6}$, by 
taking a sufficiently large $B_{5}$ including $B_{6}$, 
we may assume that $D\subseteq A\cap B_{6}'$. 
As $D$ is isomorphic to $P_{N_{1}}K(l^{2}({\bf Z}))P_{N_{1}}$, 
let $\Psi$ be the canonical isomorphism 
from $D$ onto $P_{N_{1}}K(l^{2}({\bf Z}))P_{N_{1}}$ and let
\[
f_{i}\equiv\Psi^{-1}(e_{i})
\]
for $i=0,\ldots ,m-1$. Then $(\ f_{i}\ |\ i=0,\ldots ,m-1\ )$ 
is an orthogonal family of projections in $A\cap B_{6}'$ such that
\[
\| \alpha_{1}(f_{i})-f_{i+1}\| <\varepsilon_{1}\ ,
\]
\begin{equation}
m \tau (f_{0})>(1-\varepsilon_{1})N_{1}\tau (e^{(0)})
\label{labzzz}
\end{equation}
for $i=0,\ldots ,m-1$, where $f_{m}\equiv f_{0}$. 
Thus we have for $i=0,\ldots ,m-1$,
\begin{eqnarray*}
\| \alpha_{\xi_{1}}(f_{i})-f_{i+1}\| &\leq&
\| \alpha_{\xi_{1}}(f_{i})-\alpha_{1}(f_{i})\| +
\| \alpha_{1}(f_{i})-f_{i+1}\| \\ 
&\leq& 2\| u_{1}-1\|+\varepsilon_{1}\\
&\leq& 2\cdot 4(N_{1}-1)\varepsilon_{2}+\varepsilon_{1}\ .
\end{eqnarray*}
Using the formula (\ref{hoshi}), 
the formula (\ref{labw'}) and
\begin{eqnarray*}
\| \alpha_{1}\alpha_{2}-\alpha_{2}\alpha_{1}\| &\leq&
2(\| \alpha_{1}-\alpha_{\xi_{1}}\|
+\| \alpha_{2}-\alpha_{\xi_{2}}\| )\\
&\leq&4(\| u_{1}-1\| +\| u_{2}-1\| )\ ,
\end{eqnarray*}
we can also make $\| \alpha_{\xi_{2}}(f_{i})-f_{i}\| $ 
very small.
Finally we want to estimate $\tau (f_{0})$. 
We have already three inequalities from (\ref{labz}),(\ref{labzz}) 
and (\ref{labzzz})
\[
1\leq (N_{1}N_{2}+1)\tau (p_{(0,0)})\ ,
\]
\[
\tau (e^{(0)})>(1-\varepsilon_{1})N_{2}\tau (p_{(0,0)})\ ,
\]
\[
m\tau (f_{0})>(1-\varepsilon_{1})N_{1}\tau (e^{(0)})\ .
\]
From these we obtain
\[
m\tau (f_{0})>(1-\varepsilon_{1})^{2}N_{1}N_{2}
\frac{1}{N_{1}N_{2}+1}\left( 
m\tau (f_{0})+\tau (1-\sum_{i=0}^{m-1}f_{i})
\right)\ .
\]
Since $1\ll k_{i}\ll l_{i}$, $(\, f_{i}\, 
|\, i=0,\ldots ,m-1\, )$ 
satisfies the desired conditions. \hfill \qed

\vspace{5mm}

\noindent Proof of Theorem \ref{lab1}.

Let $\alpha$ be a ${\bf Z}^{2}$-action 
on a UHF algebra $A$ which satisfies the condition (2). 
For any $m\in {\bf N}$ we take $m_{0},\ m_{1}\in {\bf N}$ 
such that $m\ll m_{1}\ll m_{0}$ and $m_{0}$ is divided by $m_{1}$. 
Furthermore for any $\varepsilon_{1} >0$ 
and finite subset $F$ of $A$, 
we take a unital full matrix subalgebra $B_{1}$ 
of $A$ such that for any $x\in F$ 
there exists $y\in B_{1}$ with $\| x-y\| <\varepsilon_{1}$. 
If we apply Lemma \ref{lab4} 
to any $n\in {\bf N}$ and $\varepsilon_{2} >0$ 
then we have an orthogonal family 
$(\ e_{i}\ |\ i=0,\ldots ,m-1\ )$ of projections in 
$A\cap B_{1}'$ satisfying 
\[
\| \alpha_{\xi_{1}}(e_{i})-e_{i+1}\| <\varepsilon_{2}\ ,
\]
\[
\| \alpha_{\xi_{2}}(e_{i})-e_{i}\| <\varepsilon_{2}\ ,
\]
\[
\tau (e_{0})\geq n\tau (1-\sum_{i=0}^{m_{0}-1}e_{i})
\]
for $i=0,\ldots ,m_{0}-1$, where $e_{m_{0}}\equiv e_{0}$. 
This is not sufficient because the sum $\sum_{i=0}^{m-1}e_{i}$ 
of the projections $(e_{i})$ may not be $1$. 
We will cope with this problem now. Put
\[
x_{1}=\sum_{i=0}^{m_{0}-1}
e_{i+1}\alpha_{\xi_{1}}(e_{i})+
\left(1-\sum_{i=0}^{m_{0}-1}e_{i}\right)
\left(1-\sum_{i=0}^{m_{0}-1}\alpha_{\xi_{1}}(e_{i})\right)\ ,
\]
\[
x_{2}=\sum_{i=0}^{m_{0}-1}
e_{i}\alpha_{\xi_{2}}(e_{i})+
\left(1-\sum_{i=0}^{m_{0}-1}e_{i}\right)
\left(1-\sum_{i=0}^{m_{0}-1}\alpha_{\xi_{2}}(e_{i})\right)
\]
and let $u_{1}|x_{1}|$ and $u_{2}|x_{2}|$ be 
the polar decompositions of $x_{1}$ and $x_{2}$ 
respectively. As in the proof of Lemma \ref{lab4} 
we can show that $u_{1}$ and $u_{2}$ are unitaries in 
$A$ satisfying
\[
\| u_{1}-1\| < 4m_{0}\varepsilon_{2}\ ,
\]
\[
Ad\, u_{1}\circ\alpha_{\xi_{1}}(e_{i})=e_{i+1}
\]
for $i=0,\ldots ,m_{0}-1$, where $e_{m_{0}}\equiv e_{0}$, 
and 
\[
\| u_{2}-1\| < 4m_{0}\varepsilon_{2}\ ,
\]
\[
Ad\, u_{2}\circ\alpha_{\xi_{2}}(e_{i})=e_{i}
\]
for $i=0,\ldots ,m_{0}-1$. Let 
$\alpha_{1}= Ad\, u_{1}\circ\alpha_{\xi_{1}}$ 
and let $\alpha_{2}= Ad\, u_{2}\circ\alpha_{\xi_{2}}$. 
Then $\alpha_{1}^{m_{0}}$ and $\alpha_{2}$ 
are automorphisms of $e_{0}Ae_{0}$. 
By Lemma \ref{lab4} there are an orthogonal family 
$(\ p_{j}\ |\ j=0,\ldots ,n-1\ )$ of projections in 
$A\cap B_{1}'$ and a positive number $c_{1}(m_{0},\varepsilon_{2})$ 
which decreases to zero as $\varepsilon_{2}\rightarrow 0$ such that
\[
\sum_{i=0}^{n-1}p_{i}\leq e_{0}\ ,
\]
\[
\| \alpha_{1}^{m_{0}}(p_{i})-p_{i+1}\| <c_{1}(m_{0},\varepsilon_{2})\ ,
\]
\[
\| \alpha_{2}(p_{i})-p_{i}\| <c_{1}(m_{0},\varepsilon_{2})
\]
for $i=0,\ldots ,n-1$, where $p_{n}\equiv p_{0}$, and
\[
\tau (p_{i})=\tau (1-e)\ ,
\]
where $e\equiv \sum_{i=0}^{m_{0}-1}e_{i}$.  
We have used the fact $\tau (e_{0})\geq n\tau (1-e)$ here. 
For any $\varepsilon_{3}>0$ and any unital full matrix subalgebra 
$B_{2}$ of $A$, by taking a sufficiently large $B_{1}$ 
and by applying \cite[Lemma\ 3.5]{kishimoto1}, 
there exists a partial isometry $v\in A\cap B_{2}'$ such that 
\[
v^{\ast}v=1-e,\ \ vv^{\ast}=p_{0}\ ,
\]
\begin{eqnarray*}
\| \alpha_{2}(v)-v\|
&\leq&\| \alpha_{2}(p_{0})-p_{0}\| +\| \alpha_{2}(1-e)-(1-e)\| +
\varepsilon_{3}\\
&\leq&c_{1}(m_{0},\varepsilon_{2})+m_{0}\varepsilon_{2}+\varepsilon_{3}\ .
\end{eqnarray*}
As before there also exists a unitary $u_{1}'\in A$ 
satisfying that $\| u_{1}'-1\| <4c_{1}(m_{0},\varepsilon_{2})$ and
\[
Ad\, u_{1}'\circ\alpha_{1}^{m_{0}}(p_{i})=p_{i+1}
\]
for $i=0,\ldots ,n-1$. Let 
$\beta= Ad\, u_{1}'\circ\alpha_{1}^{m_{0}}$ 
and let $w= n^{-\frac{1}{2}}\sum_{i=0}^{n-1}\beta^{i}(v)$. Then we have\
\[
w^{\ast}w=1-e,\ \ ww^{\ast}\leq e_{0}\ ,
\]
\begin{equation}
\| \beta (w)-w\| <n^{-\frac{1}{2}}\cdot 2\ ,
\label{3star}
\end{equation}
\[
\| \alpha_{2}(w)-w\| <c_{2}(m_{0},n,\varepsilon_{2},\varepsilon_{3})
\]
for some positive number 
$c_{2}(m_{0},n,\varepsilon_{2},\varepsilon_{3})$ 
which we can make very small. 
Furthermore by taking $B_{2}$ very large we may assume  
that $w\in A\cap B_{3}'$ for any unital full matrix subalgebra 
$B_{3}$ of $A$.
Let
\[
E_{i,j}=\left\{
\begin{array}{ll}
\alpha_{1}^{i-1}(w) \alpha_{1}^{i-2}(w) 
\cdots \alpha_{1}^{j}(w) & (if\ i>j)\\
\alpha_{1}^{i-1}(ww^{\ast}) & (if\ i=j)\\
\alpha_{1}^{i}(w)^{\ast} \alpha_{1}^{i+1}(w)^{\ast} 
\cdots \alpha_{2}^{j-1}(w)^{\ast} & (if\ i<j)
\end{array}
\right.
\]
for $0\leq i,\ j\leq m_{0}$ and let $C$ be 
the $C^{\ast}$-subalgebra of $A$ generated by 
$\{\, E_{i,j}\, |\, 0\leq i,j\leq m_{0}-1\, \}$. 
Then $C$ is isomorphic to $M_{m_{0}+1}$ 
and we may assume that $C$ is a subalgebra of $A\cap B_{4}'$ 
for any unital full matrix subalgebra $B_{4}$ of $A$ if $B_{3}$ 
is very large. Let
\[
U=\left[
\begin{array}{ccccccc}
1 &       &      &       &       &  &       \\
  &0      &\cdot &\cdots &\cdot  &0 &1      \\
  &1      &0     &\cdots &\cdot  &  &0      \\
  &0      &1     &\cdot  &       &  &\cdot  \\
  &\vdots &      &\ddots &\ddots &  &\vdots \\
  &\cdot  &      &0      &1      &0 &0      \\
  &0      &      &       &0      &1 &0      
\end{array}
\right]\in M_{m_{0}+1}\ .
\]
By simple calculation 
$\alpha_{1}\lceil C=Ad\, U$ and 
\[
{\rm Sp}(U)=\{ 1\}
\cup \{\,  e^{\frac{2\pi ik}{m_{0}}}\, |\, k=0,\ldots ,m_{0}-1\, \}\ . 
\]
Define orthogonal families 
$(\, e^{(1)}_{i}\, |\, i=0,\ldots ,m_{1}-1\, ),
\ (\, e^{(2)}_{j}\, |\, j=0,\ldots ,m_{1}\, )$ 
of projections in $A\cap B_{4}'$ as follows:
\begin{eqnarray*}
e^{(1)}_{i}
&\equiv&\sum_{k=0}^{\frac{m_{0}}{m_{1}}-2}
E_{m_{1}+i(\frac{m_{0}}{m_{1}}-1)+k,\ m_{1}+i(\frac{m_{0}}
{m_{1}}-1)+k}\ ,\\ 
e^{(2)}_{j}&\equiv&E_{j,j}\ .
\end{eqnarray*}
Since $(\, E_{i,j}\, |\, 0\leq i,j\leq m_{0}-1\, )$ 
is a system of matrix units in $A\cap B_{4}'$, 
there are canonical systems $(\, e^{(1)}_{k,l}\, 
|\, 0\leq k,l\leq m_{1}-1\, )$, $(\, e^{(2)}_{k,l}\, 
|\, 0\leq k,l\leq m_{1}\, )$ of matrix units 
associated with $(\, e^{(1)}_{k}\, )$ and $(\, e^{(2)}_{k}\, )$ 
respectively. Define a partial isometry $V$ in $A\cap B_{4}'$ by
\[
V=\sum_{i=0}^{m_{1}-1}e^{(1)}_{i+1,i}
+\sum_{j=0}^{m_{1}}e^{(2)}_{j+1,j}\ ,
\]
where $e^{(1)}_{m_{1},m_{1}-1}\equiv e^{(1)}_{0,m_{1}-1}$ 
and $e^{(2)}_{m_{1}+1,m_{1}}\equiv e^{(2)}_{0,m_{1}}$. Then
\[
\sum_{i=0}^{m_{1}-1}e^{(1)}_{i}
+\sum_{j=0}^{m_{1}}e^{(2)}_{j}=1_{C}\ \ 
\left(\ =\ \sum_{i=0}^{m_{0}}E_{i,i}\ \right)\ ,
\]
\[
Ad\, V(e^{(1)}_{i})=e^{(1)}_{i+1}\ ,
\]
\[
Ad\, V(e^{(2)}_{j})=e^{(2)}_{j+1}
\]
for $i=0,\ldots ,m_{1}-1$ and $j=0,\ldots ,m_{1}$ 
where $e^{(1)}_{m_{1}}\equiv e^{(1)}_{0}$ 
and $e^{(2)}_{m_{1}+1}\equiv e^{(2)}_{0}$. 
By a simple calculation
\[
{\rm Sp}(V)=\{\, e^{\frac{2\pi ik}{m_{0}-m_{1}}}\, 
|\, k=0,\ldots ,m_{0}-m_{1}-1\, \}\cup 
\{\, e^{\frac{2\pi il}{m_{1}+1}}\, 
|\, l=0,\ldots ,m_{1}\, \}\ .
\]
Therefore ${\rm Sp}(V)$ is very close to 
${\rm Sp}(U)$ if $m_{0}$ and $m_{1}$ are very large, 
i.e., $V$ is almost unitarily equivalent to 
$U$. Consequently we can find orthogonal families 
$(\, f^{(1)}_{i}\, |\, i=0,\ldots ,m_{1}-1\, )$, 
$(\, f^{(2)}_{j}\, |\, j=0,\ldots ,m_{1}\, )$ 
of projections in $A\cap B_{4}'$ and a small enough positive number 
$c_{3}(m_{0},m_{1},\varepsilon_{2},\varepsilon_{3})$ 
in  such a way that
\[
\sum_{i=0}^{m_{1}-1}f^{(1)}_{i}+\sum_{j=0}^{m_{1}}f^{(2)}_{j}=1_{C}\ ,
\]
\[
\| \alpha_{1}(f^{(1)}_{i})-f^{(1)}_{i+1}\| <
c_{3}(m_{0},m_{1},\varepsilon_{2},\varepsilon_{3})\ ,
\]
\[
\| \alpha_{2}(f^{(1)}_{i})-f^{(1)}_{i}\| <
c_{3}(m_{0},m_{1},\varepsilon_{2},\varepsilon_{3})\ ,
\]
\[
\| \alpha_{1}(f^{(2)}_{j})-f^{(2)}_{j+1}\| <
c_{3}(m_{0},m_{1},\varepsilon_{2},\varepsilon_{3})\ ,
\]
\[
\| \alpha_{2}(f^{(2)}_{j})-f^{(2)}_{j}\| <
c_{3}(m_{0},m_{1},\varepsilon_{2},\varepsilon_{3})
\]
for $i=0,\ldots ,m_{1}-1$ and $j=0,\ldots ,m_{1}$, 
where $f^{(1)}_{m_{1}}\equiv f^{(1)}_{0}$ and 
$f^{(2)}_{m_{1}+1}\equiv f^{(2)}_{0}$. 
Then by considering (\ref{3star}), if $n$ is very large,
\[
(\, f^{(1)}_{i}\, |\, i=0,\ldots ,m_{1}-1\, ),
\]
\[
(\, f^{(2)}_{j}\, |\, j=0,\ldots ,m_{1}\, ),
\]
\[
(\, e_{i}-\alpha_{1}^{i}(ww^{\ast})\, |\, i=0,\ldots ,m_{0}-1\, )
\]
satisfy the conditions appearing in  Remark \ref{lab-99} 
(for $n=1$). Hence $\alpha$ has the Rohlin property. \hfill \qed

\section{Conjugacy}

In this section we introduce three notions of conjugacy 
for ${\bf Z}^{N}$-actions on $C^{\ast}$-algebras 
and discuss their relationship. First we prepare some 
notations. For ${\bf Z}^{N}$-actions  $\alpha ,\beta$ 
on a unital $C^{\ast}$-algebra $A$, we write 
$\alpha \stackrel{\gamma ,\varepsilon}{\approx} \beta$ 
when
\[
\| \alpha_{\xi_{i}}-\gamma\circ\beta_{\xi_{i}}\circ\gamma^{-1}\| 
\leq\varepsilon\ \ (i=1,\ldots ,N)
\]
for $\varepsilon \geq 0$ and $\gamma\in Aut(A)$. 
For simplicity $\stackrel{\gamma ,0}{\approx}$ 
will be denoted by $\stackrel{\gamma}{\cong}$ or 
$\cong$. Recall that a 1-cocycle for $\alpha$ means 
a mapping $u$ from ${\bf Z}^{N}$ 
into the unitaries $U(A)$ of $A$ satisfying 
$u_{g+h}=u_{g}\alpha_{g}(u_{h})$ 
for each $g,h\in {\bf Z}^{N}$.

\begin{definition}\label{lab9}
{\rm 
Let $\alpha$ and $\beta$ be ${\bf Z}^{N}$-actions 
on a unital $C^{\ast}$-algebra $A$. 

\noindent (1) $\alpha$ and $\beta$ are 
{\sl approximately conjugate} 
if for any $\varepsilon >0$ there exists an automorphism 
$\gamma$ of $A$ such that 
$\alpha\stackrel{\gamma ,\varepsilon}{\approx}\beta$. 

\noindent (2) $\alpha$ and $\beta$ are 
{\sl cocycle conjugate} if there exist an automorphism 
$\gamma$ of $A$ and a 1-cocycle $u$ for $\alpha$ such that
\[
Ad\, u_{g}\circ\alpha_{g}=\gamma\circ\beta_{g}\circ\gamma^{-1}
\]
for each $g\in {\bf Z}^{N}$.

\noindent (3) $\alpha$ and $\beta$ are 
{\sl outer conjugate} 
if there exist an automorphism $\gamma$ of $A$ 
and unitaries  $u_{1},\ldots ,u_{N}$ in $A$ such that
\[
Ad\, u_{i}\circ\alpha_{\xi_{i}}
=\gamma\circ\beta_{\xi_{i}}\circ\gamma^{-1}
\]
for $i=1,\ldots ,N$ .
}
\end{definition}

Cocycle conjugacy of course implies outer conjugacy, and we have

\begin{proposition}\label{cocyc5}
Assume that $A$ is a simple separable unital $C^{\ast}$-algebra 
with a unique tracial state. Then approximately conjugate 
${\bf Z}^{N}$-actions 
on $A$ are cocycle conjugate.
\end{proposition}

Our proof is based on the generalization 
of the determinant introduced by P. de la Harpe and G. Skandalis 
(see \cite{hs} for details) and the famous $2\times 2$ 
matrix trick due to A. Connes (see \cite{connes0}). 
We quickly review basic facts on the former. 
For a unital $C^{\ast}$-algebra $A$, 
we let $GL_{n}(A)$ the group of 
the invertible elements in the $n\times n$ matrices $M_{n}(A)$ 
over $A$ (equipped with the $C^{\ast}$-norm). 
The inductive limit of topological groups 
$(\, GL_{n}(A)\, |\, n\in {\bf N}\, )$ 
with the usual embeddings 
$GL_{n}(A)\hookrightarrow GL_{n+1}(A)$ is denoted by 
$GL_{\infty}(A)$ and the connected component 
of the identity by $GL_{\infty}^{0}(A)$. 
Suppose that $\tau$ is a tracial state on 
$A$. If $\xi$ is a piecewise continuously differentiable mapping 
from $[0,1]$ into $GL_{\infty}^{0}(A)$, we define
\[
{\widetilde{\Delta}}_{\tau}(\xi )=
\frac{1}{2\pi i}\int_{0}^{1}
\tau\left( \dot{\mbox{$\xi$}}(t)\xi (t)^{-1}\right)dt
\]
(note that the range of $\xi$ is contained in $GL_{n}(A)$ 
for some $n$ since [0,1] is compact and that $\tau$ 
actually means $\tau\otimes {\rm tr}$ 
on $A\otimes M_{n}=M_{n}(A)$). 
The determinant $\Delta_{\tau}$ (\cite{hs}) associated 
with a tracial state $\tau$ is the mapping 
from $GL_{\infty}^{0}(A)$ into ${\bf C}/\tau_{\ast}(K_{0}(A))$ 
defined by
\[
\Delta_{\tau}(x)=p({\widetilde{\Delta}}_{\tau}(\xi ))\ .
\]
Here $p$ is the quotient mapping 
from ${\bf C}$ onto ${\bf C}/\tau_{\ast}(K_{0}(A))$, 
and $\xi$ is a piecewise continuously differentiable mapping 
from [0,1] into $GL_{\infty}^{0}(A)$ with $\xi (0)=1$ 
and $\xi (1)=x$.

A crucial fact here is that $\Delta_{\tau}$ 
is a group homomorphism. For a unitary $x\in A$ with $\| x-1\| <1$, 
the logarithm $h=i^{-1}\log (x)$\,(with the principal branch) 
makes sense and we can consider the path $\xi(t)=\exp (iht)$ 
from $1$ to $x$. 
Since $\dot{\mbox{$\xi$}}(t)\xi (t)^{-1}=ih$, we have 
\[
\Delta_{\tau}(x)=
\frac{1}{2\pi i}p\left( \int_{0}^{1}\tau (ih)dt\right)=
\frac{1}{2\pi i}p(\tau (\log (x)))\ .
\]

\vspace{4mm}

\noindent Proof of Proposition \ref{cocyc5}

Let $\alpha$ and $\beta$ be approximately conjugate 
${\bf Z}^{N}$-actions 
on $A$. 
In general, if an automorphism 
of a simple unital $C^{\ast}$-algebra is close 
to the identity in norm 
then it is inner, and furthermore 
it is implemented by a unitary which is also close to the unit 
of the algebra. Hence for a sufficiently small $\varepsilon >0$ 
there exist unitaries $u_{1},\ldots ,u_{N}$ in $A$ 
and an automorphism $\gamma$ of $A$ such that
\[
Ad\, u_{i}\circ\alpha_{\xi_{i}}
=\gamma\circ\beta_{\xi_{i}}\circ\gamma^{-1}\ ,
\]
\[
\| u_{i}-1\| <\varepsilon
\]
for $i=1,\ldots ,N$. We want to show 
\begin{equation}
u_{k}\alpha_{\xi_{k}}(u_{l})=u_{l}\alpha_{\xi_{l}}(u_{k})
\label{cocyc2}
\end{equation}
for any $k,l=1,\ldots ,N$. From the commutativity 
of $\alpha_{\xi_{k}}$, 
and the simplicity of $A$, there exists $\lambda\in{\bf T}$ 
such that
\[
(u_{l}\alpha_{\xi_{l}}(u_{k}))^{\ast}u_{k}\alpha_{\xi_{k}}(u_{l})
=\lambda 1\ .
\]
Since $u_{k},u_{l}$ are close to $1$, we can set
\[
h_{k}=\frac{1}{2\pi i}\log (u_{k}),\ \ h_{l}
=\frac{1}{2\pi i}\log (u_{l})
\]
and
\[
H(s)=\frac{1}{2\pi i}\log\left\{ 
(u_{l}^{s}\alpha_{\xi_{l}}(u_{k}^{s}))^{\ast}
u_{k}^{s}\alpha_{\xi_{k}}(u_{l}^{s})
\right\}
\]
for $s\in [0,1]$. Applying $\Delta_{\tau}$ 
to the both sides of the equality
\[
e^{-2\pi i\alpha_{\xi_{l}}(sh_{k})}e^{-2\pi ish_{l}}
e^{2\pi ish_{k}}e^{2\pi i\alpha_{\xi_{k}}(sh_{l})}=
e^{2\pi iH(s)}\ ,
\]
we have
\[
-p(\tau (\alpha_{\xi_{l}}(sh_{k})))-p(\tau (sh_{l}))
+p(\tau (sh_{k}))+p(\tau (\alpha_{\xi_{k}}(sh_{l})))=
p(\tau (H(s)))\ . 
\]
The uniqueness of a trace shows that the left hand side 
of the above equality is zero, i.e., 
$\tau (H(s))\in \tau_{\ast}(K_{0}(A))$ 
for any $s\in [0,1]$. Since $\tau_{\ast}(K_{0}(A))$ 
is discrete in ${\bf C}$ and $\tau (H(0))=0$, 
we obtain $H(1)=\tau (H(1))=0$\ \ i.e., $\lambda =1$. 
Using these unitaries $u_{1},\ldots ,u_{N}$, 
we construct a desired 1-cocycle by the method 
in \cite{connes0}. 
We consider the ${\bf Z}^{N}$-action $\sigma$ on $M_{2}(A)$ 
defined by 
\[
\sigma_{\xi_{i}}=Ad\, 
\left[
\begin{array}{cc}
1 & 0 \\
0 & u_{i}
\end{array}
\right]
\circ\alpha_{\xi_{i}}\ .
\]
Since $\sigma_{\xi_{1}},\ldots ,\sigma_{\xi_{N}}$ 
commute with each others from (\ref{cocyc2}), 
$\sigma$ is indeed well-defined. Note that
\[
\sigma_{g}
\left(\left[
\begin{array}{cc}
x & 0 \\
0 & y
\end{array}
\right]\right)
=
\left[
\begin{array}{cc}
\alpha_{g}(x) & 0 \\
0 & \gamma\circ\beta_{g}\circ\gamma^{-1}(y)
\end{array}
\right]
\]
for any $g\in{\bf Z}^{N}$ and $x,y\in A$. The identity
\[
\left[
\begin{array}{cc}
0 & 0 \\
0 & x
\end{array}
\right]
=
\left[
\begin{array}{cc}
0 & 0 \\
1 & 0
\end{array}
\right]
\left[
\begin{array}{cc}
x & 0 \\
0 & 0
\end{array}
\right]
\left[
\begin{array}{cc}
0 & 1 \\
0 & 0
\end{array}
\right]
\]
shows 
\[
\gamma\circ\beta_{g}\circ\gamma^{-1}(x)
=Ad\, u_{g}\circ\alpha_{g}(x)\ ,
\]
where $u_{g}$ is the desired 1-cocycle defined by
\[
\left[
\begin{array}{cc}
0 & 0 \\
u_{g} & 0
\end{array}
\right]
=
\sigma_{g}\left(\left[
\begin{array}{cc}
0 & 0 \\
1 & 0
\end{array}
\right]\right)
\]
(see \cite[Lemma\ 8.11.2]{pedersen}) .  \hfill\qed

\begin{remark}\label{lab9.1}
{\rm 
If $\alpha$ and $\beta$ are outer conjugate automorphisms 
of a UHF algebra with the Rohlin property, 
then they are approximately conjugate 
by the stability property. 
Hence the three notions of conjugacy defined above 
are the same for those automorphisms. 
However outer conjugacy does not imply 
approximate conjugacy for ${\bf Z}^{N}$-actions. 
See Remark \ref{lab13} for a counter-example.
}
\end{remark}

\section{Product type actions}

In this section we discuss product type ${\bf Z}^{2}$-actions 
on UHF algebras. As in the case 
of single automorphisms, the Rohlin property 
for these actions is closely related to a notion 
of uniform distribution of points in ${\bf T}^{2}$. 
First we say the $N$-dimensional version of 
\cite[Lemma\ 4.1]{bek}. This can be shown as in the one-dimensional 
case, so we omit the proof.

\begin{proposition}\label{lab5}\ \ 
Let $(\, S_{k}\, |\, k\in {\bf N}\, )$ be a sequence 
of finite sequences in ${\bf T}^{N}$ i.e., 
\[
S_{k}=(\, s_{k,p}\, |\, p=1,\ldots ,n_{k}\, )\ ,
\]
\[
s_{k,p}\in {\bf T}^{N}
\]
for each $k\in {\bf N}$ and $p=1,\ldots ,n_{k}$. 
Then the following conditions on 
$(\, S_{k}\, |\, k\in{\bf N}\, )$ 
are equivalent.

\vspace{2mm}

\noindent {\rm (1)}
\[
\lim_{k\rightarrow\infty}\frac{1}{n_{k}}
\sum_{p=1}^{n_{k}}f(s_{k,p})=
\int_{{\bf T}^{n}}f(s)ds
\]
for any $f\in C({\bf T}^{N})$, where $ds$ 
denotes the normalized Haar measure on ${\bf T}^{N}$.

\noindent {\rm (2)}
\[
\lim_{k\rightarrow\infty}\frac{1}{n_{k}}
\sum_{p=1}^{n_{k}}s_{k,p}^{l}=0
\]
for any $l=(l_{1},\ldots ,l_{N})\in {\bf Z}^{N}
\setminus\{0\}$, where $s^{l}$ denotes $s_{1}^{l_{1}}
\cdots s_{N}^{l_{N}}$ for each $s=(s_{1},\ldots ,s_{N})
\in {\bf T}^{N}$.

\noindent {\rm (3)}
\[
\lim_{k\rightarrow\infty}\frac{1}{n_{k}}\nu_{k}
\left(\prod_{i=1}^{N}[\,\theta_{1}^{(i)},
\theta_{2}^{(i)})\right)=
(2\pi)^{-N}\prod_{i=1}^{N}(\theta_{2}^{(i)}
-\theta_{1}^{(i)})
\]
for any $0\leq\theta_{1}^{(i)}\leq\theta_{2}^{(i)}<2\pi$, 
where $\nu_{k}$ is defined by
\[
\nu_{k}(S)= \sharp \{\ p\ |\ 1\leq p\leq n_{k}\ 
and\ {\rm arg}(s_{k,p})\in S\ \}
\]
for each subset $S$ of $\prod_{i=1}^{N}[\, 0,2\pi )$ 
and $\sharp F$ denotes the cardinality of the set $F$. 

\vspace{4mm}

\noindent These conditions necessarily imply that 
$n_{k}\rightarrow\infty$.
Moreover suppose that 
$(\, n_{k}\, |\, k\in {\bf N}\,)$ 
has the following 
asymptotic factorization into large factors: 
\ For any $n\in{\bf N}$ 
there exists positive integer $k_{0}$ such that 
for any positive integer $k\geq k_{0}$ 
one has 
$n^{(1)}_{k},\ldots ,n^{(N)}_{k}\geq n$ 
where $n^{(i)}_{k}$ are the components of $n_{k}$, i.e., 
$n_{k}=n^{(1)}_{k}\cdots n^{(N)}_{k}$. 
Then the above conditions are also equivalent to

\vspace{2mm}

\noindent {\rm (4)} For any $\varepsilon >0$ 
there exist positive integers $k_{0}$ and $n_{0}$ 
such that for any 
$k,\ n_{k}^{(1)},\ldots ,n_{k}^{(N)}\in {\bf N}$ 
satisfying 
$k\geq k_{0},\ n_{k}^{(1)},\ldots ,n_{k}^{(N)}\geq n_{0}$ 
and $n_{k}=n_{k}^{(1)}\cdots n_{k}^{(N)}$, 
there exists an bijection $\varphi$ 
from $\{ 1,\ldots ,n_{k}\}$ 
onto 
$\{ 1,\ldots ,n_{k}^{(1)}\}\times\cdots\times 
\{ 1,\cdots ,n_{k}^{(N)}\}$ 
such that
\begin{equation}
\left|\, 
s_{k,p}-
({\rm exp}\left( 2\pi i\cdot\frac{(\varphi (p))_{1}}
{n_{k}^{(1)}}\right),\ldots ,
{\rm exp}\left( 2\pi i\cdot\frac{(\varphi (p))_{N}}
{n_{k}^{(N)}}\right))
\right|\, <\varepsilon
\label{labdis}
\end{equation}
for any $k$ and $p$, where $|s|\equiv {\rm max}
\{\, |s_{p}|\,:1\leq p\leq N\,\}$ for each 
$s\in {\bf T}^{N}$ and $(\varphi (p))_{i}$ 
denotes the i-th component of $\varphi (p)$. 
\end{proposition}

If $S_{k}$ satisfies the estimate (\ref{labdis}) 
for some $\varphi$ as above, then $S_{k}$ is said 
to be $(n_{k}^{(1)},\ldots ,n_{k}^{(N)};\varepsilon)$
-{\sl distributed}. If one of the conditions 
of the above proposition holds then 
$(\, S_{k}\, |\, k\in {\bf N}\, )$ 
is said to be {\sl uniformly distributed}.

\begin{definition}\label{lab6}
{\rm 
Let $\alpha$ be a ${\bf Z}^{N}$-action 
on a UHF algebra $A$. Then $\alpha$ is said 
to be a {\sl product type action} 
if there exists a sequence $(\, m_{k}\, |\, k\in {\bf N}\, )$ 
of positive integers such that 
$A\cong \otimes_{k=1}^{\infty}M_{m_{k}}$ and
\[
\alpha_{g}(A_{k})=A_{k}
\]
for any $g\in {\bf Z}^{N}$ and $k\in {\bf N}$, 
where $A_{k}$ denotes the $C^{\ast}$-subalgebra 
of $A$ corresponding to 
$M_{m_{k}}\otimes (\otimes_{l\neq k}{\bf C}1_{m_{l}})$.
}
\end{definition}

\begin{remark}\label{lab30}
{\rm 
In the situation above, if $N=2$, 
then one finds unitaries $u^{(1)}_{k},u^{(2)}_{k}$ 
in $A_{k}$ and $\lambda_{k}\in {\bf T}$ such that
\[
\alpha_{(p,q)}\lceil A_{k}=Ad\, {u^{(1)}_{k}}^{p}{u^{(2)}_{k}}^{q}\ ,
\]
\[
u^{(1)}_{k}u^{(2)}_{k}=\lambda_{k}u^{(2)}_{k}u^{(1)}_{k}
\]
for any $p,q\in{\bf Z}$. Since $u^{(1)}_{k},u^{(2)}_{k}$ are unique 
up to a constant multiple, $\lambda_{k}$ is unique. 
In addition $\lambda_{k}^{m_{k}}=1$. For if $\mu_{1}$ 
is an eigenvalue of $u^{(2)}_{k}$ with multiplicity $r_{1}$ 
then $\mu_{1}\lambda_{k}^{-p}$ is also an eigenvalue 
of $u^{(2)}_{k}$ with multiplicity $r_{1}$ for each 
$p\in {\bf N}$. Since $M_{m_{k}}$ is finite-dimensional, 
there exists $p_{0}\in {\bf N}$ such that 
$\lambda_{k}^{p_{0}}=1$ and $\lambda_{k}^{p}\neq 1$ 
for any $p=1,\ldots ,p_{0}-1$. If 
$\{\, \mu_{1}\lambda_{k}^{-p}\, |\, p=0,\ldots ,p_{0}-1\, \}$ 
does not exhaust all the eigenvalues of $u^{(2)}_{k}$ 
then we take an eigenvalue $\mu_{2}$ of $u^{(2)}_{k}$ 
not  belonging to 
$\{\, \mu_{1}\lambda_{k}^{-p}\, |\, p=0,\ldots ,p_{0}-1\, \}$ 
and repeat the same process. Thus there exist eigenvalues 
$\mu_{1},\ldots ,\mu_{s}$ of $u^{(2)}_{k}$ with multiplicity 
$r_{1},\ldots ,r_{s}$ respectively. 
Since $m_{k}=(r_{1}+\cdots +r_{s})p_{0}$, 
it follows that $\lambda_{k}^{m_{k}}=1$.
}
\end{remark}

For  $n\times n$ unitary matrices $U$ and $V$ with $UV=VU$, 
we define ${\rm Sp}(U)$ to be a sequence consisting 
of the eigenvalues of $U$, each repeated as often 
as multiplicity dictates and ${\rm Sp}(U,V)$ is 
a sequence consisting of the pairs of eigenvalues 
of $U$ and $V$ with a common eigenvector, each repeated 
as often as multiplicity dictates. 
Then the Rohlin property for the product type 
${\bf Z}^{2}$-actions on $A$ with $\lambda_{k}=1$ 
can be characterized as follows.

\begin{proposition}\label{lab8}
Let $\alpha$ be a product type ${\bf Z}^{2}$-action 
on a UHF algebra $A$ with 
$(\, m_{k}\, |\, k\in {\bf N}\, )$, 
$(\, u^{(1)}_{k}\, |\, k\in {\bf N}\, )$, 
$(\, u^{(2)}_{k}\, |\, k\in {\bf N}\, )$, 
$(\,\lambda_{k}\, |\, k\in {\bf N}\, )$ as above. 
If $\lambda_{k}=1$ for each $k\in {\bf N}$ 
then the following conditions are equivalent:
\begin{list}{}{}
\item[{\rm (1)}] $\alpha$ has the Rohlin property.
\item[{\rm (2)}] $(\ {\rm Sp}(\otimes_{k=m}^{n}u^{(1)}_{k},
\otimes_{k=m}^{n}u^{(2)}_{k})\ |\ n=m,\ m+1,\ \ldots\ )$ 
is uniformly distributed for any $m\in {\bf N}$\ .
\end{list}
\end{proposition}

\noindent {\it Proof}.\ \ By Corollary \ref{lab23}, 
(1) is equivalent to the condition: 
$\alpha_{\xi_{1}}^{p}\alpha_{\xi_{2}}^{q}$ 
has the Rohlin property as a single automorphism 
of $A$ for each $(p,q)\in{\bf Z}^{2}\setminus\{ 0\}$. 
By \cite[Lemma\ 5.2]{kishimoto2} 
this condition is equivalent to the condition: 
\[
(\ {\rm Sp}
\left(\bigotimes_{k=m}^{n}{u^{(1)}_{k}}^{p}{u^{(2)}_{k}}^{q}\right)\ 
|\ n=m,\ m+1,\ \ldots\ )
\]
is uniformly distributed in ${\bf T}$.  
By Proposition \ref{lab5} for $N=1$, 
the last condition is equivalent to the condition: 
for any $m\in {\bf N}$
\[
\lim_{n\rightarrow\infty}\frac{1}{N(m,n)}
\sum_{(\lambda_{1},\lambda_{2})\in {\rm Sp}
\left(\otimes_{k=m}^{n}u^{(1)}_{k},\otimes_{k=m}^{n}u^{(2)}_{k}\right)}
\lambda_{1}^{p}\lambda_{2}^{q}=0\ ,
\]
where $N(m,n)\equiv\prod_{k=m}^{n}m_{k}$. 
Finally by Proposition \ref{lab5} for $N=2$, 
the last condition is equivalent to (2). \hfill \qed

\vspace{5mm}

In \cite{kishimoto2} A.Kishimoto showed the following 
for a UHF algebra $A$.
\begin{list}{}{}
\item[{\rm (1)}] Any two product type ${\bf Z}$-actions 
on $A$ with the Rohlin property  
are approximately conjugate.
\item[{\rm (2)}] For any ${\bf Z}$-action $\alpha$  
on $A$ with the Rohlin property and $\varepsilon >0$, 
there exist a product type ${\bf Z}$-action $\beta$ 
on $A$ with the Rohlin property and an automorphism 
$\gamma$ of $A$ such that
$\alpha
\stackrel{\gamma ,\varepsilon}{\approx}
\beta$.
\end{list}
In particular there is one and only one 
approximate conjugacy class of ${\bf Z}$-actions
on $A$ with the Rohlin property. In the case of $N=2$ 
we do not know whether (2) is valid or not. 
In the rest of this section we state several versions of (1) 
for ${\bf Z}^{2}$.

\begin{theorem}\label{lab10}
Let $\alpha$ and $\beta$ be product type ${\bf Z}^{2}$-actions  
on a UHF algebra $A$ 
with the Rohlin property. Let $\alpha$ be determined by 
$(\, m_{k}\, |\, k\in{\bf N}\, )$, 
$(\, \lambda_{k}\, |\, k\in{\bf N}\, )$ 
as in Definition \ref{lab6} and Remark \ref{lab30}, 
and $\beta$ by 
$(\, n_{l}\, |\, l\in{\bf N}\, )$, 
$(\, \mu_{l}\, |\, l\in{\bf N}\, )$. 
If $\lambda_{k}=\mu_{l}=1$ for each $k,l\in {\bf N}$, 
then $\alpha$ and $\beta$ are approximately conjugate.\end{theorem}

\noindent {\it Proof}.\ \ By patching several parts 
of the $M_{m_{k}}$'s and the $M_{n_{l}}$'s 
respectively there exist a sequence $(\, N_{k}\, |\, k\in {\bf N}\, )$ 
of integers satisfying $N_{0}=0$ and $N_{k}>0$ ($k\geq 1$) 
and sequences $(\, U^{(1)}_{k}\, |\, k\in {\bf N}\, )$, 
$(\, U^{(2)}_{k}\, |\, k\in {\bf N}\, )$ of unitary matrices such that
\[
U^{(i)}_{k}\in U(M_{N_{k}}\otimes M_{N_{k+1}})\ ,
\]
\[
U^{(1)}_{k}U^{(2)}_{k}=U^{(2)}_{k}U^{(1)}_{k}\ \ 
(i=1,2,\ k\in {\bf N})\ ,
\]
\[
(A,\alpha_{\xi_{1}},\alpha_{\xi_{2}})\cong
\left(\bigotimes_{k=0}^{\infty}M_{N_{k}},
\bigotimes_{k=0}^{\infty}Ad\, U^{(1)}_{2k},
\bigotimes_{k=0}^{\infty}Ad\, U^{(2)}_{2k}\right)\ ,
\]
\[
(A,\beta_{\xi_{1}},\beta_{\xi_{2}})\cong
\left(\bigotimes_{k=0}^{\infty}M_{N_{k}},
\bigotimes_{k=0}^{\infty}Ad\, U^{(1)}_{2k+1},
\bigotimes_{k=0}^{\infty}Ad\, U^{(2)}_{2k+1}\right)\ .
\]
Let $\varepsilon >0$. Since 
$(\ {\rm Sp}\, \left(\, \otimes_{k=0}^{n}
Ad\, U^{(1)}_{2k+1},\otimes_{k=0}^{n}Ad\, U^{(2)}_{2k+1}\right)\, 
|\, n\in {\bf N}\, )$ is uniformly distributed 
by Proposition \ref{lab8} and 
$U^{(1)}_{0}U^{(2)}_{0}=U^{(2)}_{0}U^{(1)}_{0}$, 
for a sufficiently large $n$ there exist a unitary 
$W_{1}$ in $\otimes_{k=0}^{2n+2}M_{N_{k}}$ 
and unitaries $V^{(1)}_{2},V^{(2)}_{2}$ in 
$\otimes_{k=2}^{2n+2}M_{N_{k}}$ such that
\[
\|W_{1}\left(\bigotimes_{k=0}^{n}U^{(1)}_{2k+1}\right) W_{1}^{\ast}-
U^{(1)}_{0}\otimes V^{(1)}_{2}\| <2^{-2}\varepsilon\ ,
\]
\[
\|W_{1}\left(\bigotimes_{k=0}^{n}U^{(2)}_{2k+1}\right) W_{1}^{\ast}-
U^{(2)}_{0}\otimes V^{(2)}_{2}\| <2^{-2}\varepsilon\ .
\]
Then replace 
$M_{N_{2}}$ by $\bigotimes_{k=2}^{2n+2}M_{N_{k}}$
and 
$U^{(i)}_{1}$ by $\bigotimes_{k=0}^{n}U^{(i)}_{2k+1}$, 
$U^{(i)}_{2}$ by $\bigotimes_{k=1}^{n+1}U^{(i)}_{2k}$ 
($i=1,2$) respectively, and further replace 
$M_{N_{k-2n}}$ by $M_{N_{k}}$ 
($k\geq 2n+3$) and 
$U^{(i)}_{2(k-n)+1}$ by $U^{(i)}_{2k+1}$, 
$U^{(i)}_{2(k-n)+2}$ by $U^{(i)}_{2k+2}$ 
($i=1,2,\ k\geq n+1$) respectively. Thereby 
$W_{1}\in U(M_{N_{1}}\otimes M_{N_{2}})$ 
and $V^{(1)}_{2},V^{(2)}_{2}\in U(M_{N_{2}})$, 
which satisfy 
\[
\| Ad\, W_{1}(U^{(1)}_{1})-U^{(1)}_{0}
\otimes V^{(1)}_{2}\| <2^{-2}\varepsilon\ ,
\]
\[
\| Ad\, W_{1}(U^{(2)}_{1})-U^{(2)}_{0}
\otimes V^{(2)}_{2}\| <2^{-2}\varepsilon\ .
\]
In the same way, after replacing $M_{N_{3}}
,\ U^{(i)}_{2},\ U^{(i)}_{3}$ etc. suitably, 
there exist a unitary 
$W_{2}$ in $M_{N_{2}}\otimes M_{N_{3}}$ 
and unitaries $V^{(1)}_{3},V^{(2)}_{3}$ in $M_{N_{3}}$ 
such that
\[
\| Ad\, W_{2}(U^{(1)}_{2})-V^{(1)}_{2}
\otimes V^{(1)}_{3}\| <2^{-3}\varepsilon\ ,
\]
\[
\| Ad\, W_{2}(U^{(2)}_{2})-V^{(2)}_{2}
\otimes V^{(2)}_{3}\| <2^{-3}\varepsilon\ .
\]
By repeating the above procedure for $k=3,4,\ldots$, 
we can construct a unitary $W_{k}$ in 
$M_{N_{k}}\otimes M_{N_{k+1}}$ 
and unitaries $V^{(1)}_{k+1},V^{(2)}_{k+1}$ in 
$M_{N_{k+1}}$ in such a way that
\[
\| Ad\, W_{k}(U^{(1)}_{k})-V^{(1)}_{k}
\otimes V^{(1)}_{k+1}\| <2^{-(k+1)}\varepsilon\ ,
\]
\[
\| Ad\, W_{k}(U^{(2)}_{k})-V^{(2)}_{k}
\otimes V^{(2)}_{k+1}\| <2^{-(k+1)}\varepsilon\ .
\]
Thus 
\begin{eqnarray*}
(A,\alpha_{\xi_{1}},\alpha_{\xi_{2}})
&\cong &
\left(\ \bigotimes_{k=0}^{\infty}M_{N_{k}},\ 
\bigotimes_{k=0}^{\infty}Ad\, U^{(1)}_{2k},\ 
\bigotimes_{k=0}^{\infty}Ad\, U^{(2)}_{2k}\ \right) \\
 &\stackrel{\gamma_{e},\frac{1}{3}\varepsilon }{\approx}&
\left(\ \bigotimes_{k=0}^{\infty}M_{N_{k}},\ 
\bigotimes_{k=0}^{\infty}Ad\, V^{(1)}_{k},\ 
\bigotimes_{k=0}^{\infty}Ad\, V^{(2)}_{k}\ \right)\ ,
\end{eqnarray*}
where $\gamma_{e}\equiv\otimes_{k=1}^{\infty}Ad\, W_{2k}$ and

\begin{eqnarray*}
(A,\beta_{\xi_{1}},\beta_{\xi_{2}})
&\cong &
\left(\ \bigotimes_{k=0}^{\infty}M_{N_{k}},\ 
\bigotimes_{k=0}^{\infty}Ad\, U^{(1)}_{2k+1},\ 
\bigotimes_{k=0}^{\infty}Ad\, U^{(2)}_{2k+1}\ \right)\\
&\stackrel{\gamma_{o},\frac{2}{3}\varepsilon }{\approx}&
\left(\ \bigotimes_{k=0}^{\infty}M_{N_{k}},\ 
\bigotimes_{k=0}^{\infty}Ad\, V^{(1)}_{k},\ 
\bigotimes_{k=0}^{\infty}Ad\, V^{(2)}_{k}\ \right)\ ,
\end{eqnarray*}
where $\gamma_{o}\equiv\otimes_{k=0}^{\infty}Ad\, W_{2k+1}$. 
This completes the proof. \hfill \qed

\vspace{4mm}

As mentioned in the introduction, 
we discuss product type actions for two classes 
of UHF algebras. Let $(\, p_{k}\, |\, k\in {\bf N}\, )$ 
be the prime numbers in the increasing order. 
For a sequence $(\, i_{k}\, |\, k\in {\bf N}\, )$ 
of nonnegative integers with 
$\sum_{k=1}^{\infty}i_{k}=\infty$, put $q_{k}= p_{k}^{i_{k}}$ 
and let $A= \otimes_{k=1}^{\infty}M_{q_{k}}$. 
We regard $M_{q_{k}}$ as a $C^{\ast}$-subalgebra 
of $A$. We consider the class of product type ${\bf Z}^{2}$-actions  
$\alpha$ on this $A$. 
Assume that $\alpha$ looks like 
\begin{equation}
\alpha_{(p,q)}\lceil M_{q_{k}}=
Ad\, {u^{(1)}_{k}}^{p}{u^{(2)}_{k}}^{q} 
\label{5star}
\end{equation}
on $M_{q_{k}}$ with unitaries 
$u^{(1)}_{k},u^{(2)}_{k}$ in $M_{q_{k}}$ 
and $\lambda_{k}\in {\bf T}$ satisfying 
$u^{(1)}_{k}u^{(2)}_{k}=\lambda_{k}u^{(2)}_{k}u^{(1)}_{k}$. 
Since $\lambda_{k}^{q_{k}}=1$, we may regard 
$\lambda_{k}$ as an element of $G_{k}\equiv {\bf Z}/q_{k}{\bf Z}$. 
We let $[\alpha ]$ be the sequence 
$(\, \lambda_{k}\, |\, k\in {\bf N}\, )$ 
in 
$\prod_{k=1}^{\infty}G_{k}$. We define 
an equivalence relation in 
$\prod_{k=1}^{\infty}G_{k}$ by: 
$g\sim h$ if there is an $n$ such that 
$g_{k}=h_{k}$ for all $k\geq n$. Let 
$0$ be the trivial sequence $(\, 0,0,\ldots\, )$. 
We note that for every $g\in \prod_{k=1}^{\infty}G_{k}$ 
there is an action $\alpha$ in the above class with 
$[\alpha ]=g$.

\begin{theorem}\label{lab10.2}
{\rm (1)} If $\alpha$ is an action in the above class and 
$[\alpha ]\not\sim 0$, then $\alpha$ has the Rohlin property.\\
\noindent {\rm (2)} If $\alpha$ and $\beta$ 
are actions in the above class and satisfy 
the Rohlin property, then the following are equivalent:
\begin{list}{}{}
\item[{\rm (2.1)}] $[\alpha ]\sim [\beta]$.
\item[{\rm (2.2)}] $\alpha$ and $\beta$ are 
outer conjugate.
\end{list}
\end{theorem}

Before proving Theorem \ref{lab10.2} we introduce some notations 
and prepare a lemma. For a positive integer $n$ and 
$\lambda\in{\bf T}$ with $\lambda^{n}=1$, 
we define the $n\times n$ unitary matrices 
$S(n)$ and $\Omega (n,\lambda)$ by
\[
S(n)=\left[
\begin{array}{cccc}
0 & \cdots  & 0 & 1 \\
1 & \ddots & & 0 \\
 & \ddots & \ddots & \vdots \\
 & & 1 & 0
\end{array}
\right]\ ,\ \ 
\Omega (n,\lambda )=\left[
\begin{array}{cccc}
1 & & & \\
 & \lambda & & \\
 & & \ddots & \\
 & & & \lambda^{n-1}
\end{array}
\right]\ .
\]

\begin{lemma}\label{lab11}
Suppose that $U$ and $V$ are $n\times n$ unitary matrices such that
\[
UV={\rm exp}\left( 2\pi ik/n\right)VU
\]
for some $k\in{\bf N}$. 
Let $q/p$ be the irreducible form of $k/n$. 
Then there exist $\omega_{i} ,\mu_{i}\in{\bf T}$  
$(i=1,\ldots ,n/p)$ such that $(U,V)$ is conjugate to 
$(U_{1}\otimes U_{2},V_{1}\otimes V_{2})$, where
\[
\lambda= {\rm exp}(2\pi ik/n )
\]
\[
U_{1}=S(p)\ ,\ \ V_{1}= \Omega (p,\lambda^{-1})\ ,
\]
\[
U_{2}=\bigoplus_{i=1}^{n/p}\omega_{i}\ ,\ \ 
V_{2}=\bigoplus_{i=1}^{n/p}\mu_{i}\ .
\]
Moreover each $\omega_{i},\mu_{i}$ are 
unique up to multiples of powers of $\lambda$.
\end{lemma}

\noindent {\it Proof.}\ \ 
Since $U^{p}V=VU^{p}$ we have a complete orthonomal system 
of ${\bf C}^{n}$ consisting of the common eigenvectors of 
$U^{p}$ and $V$. We take such a system 
$(\xi_{(\kappa ,\mu)}|\kappa ,\mu)$ i.e.,
\[
U^{p}\xi_{(\kappa ,\mu)}=\kappa\xi_{(\kappa ,\mu)}\ ,
\]
\[
V\xi_{(\kappa ,\mu)}=\mu\xi_{(\kappa ,\mu)}.
\]
Then if $\omega^{p}=\kappa$, the space spanned by
\[
\xi_{(\kappa ,\mu )},\ \ \omega U\xi_{(\kappa ,\mu )},\ \ \ldots ,\ \ 
\omega^{p-1}U^{p-1}\xi_{(\kappa ,\mu )}
\]
is invariant under $U,V$, and the matrix representation 
of $(U,V)$ with respect to the above basis is 
$(\omega U_{1},\mu V_{1})$. Thus $(U,V)$ 
is conjugate to the direct sum of 
$(\omega_{i}U_{1},\mu_{i}V_{1})$ for some sequences 
$\omega_{i},\mu_{i}$ in ${\bf T}$. Since 
$(\omega U_{1},\mu V_{1})$ is conjugate to 
$(\omega\lambda^{k} U_{1},\mu\lambda^{j} V_{1})$ 
for all $k,j$, the last statement is obvious. \hfill \qed

\vspace{4mm}

\noindent Proof of Theorem \ref{lab10.2} (1)

Let $\alpha$ be given as in the theorem. 
Take unitaries $u^{(1)}_{k},u^{(2)}_{k}$ in $M_{q_{k}}$ 
and $\lambda_{k}\in {\bf T}$ as in (\ref{5star}). 
By Corollary \ref{lab23} it suffices to prove that 
$\alpha_{(p,q)}$ has the Rohlin property 
as a single automorphism for each $(p,q)\in {\bf Z}^{2}
\setminus \{ 0\}$. 
From the assumption there is a subsequence 
$(\, p_{k_{n}}\, |\, n\in {\bf N}\, )$ of 
$(\, p_{k}\, |\, k\in {\bf N}\, )$ such that 
$\lambda_{k_{n}}\neq 1$ for any $n$. Applying Lemma \ref{lab11} 
to $u^{(1)}_{k_{n}},u^{(2)}_{k_{n}}$ we have a decomposition 
$(u^{(1)}_{k_{n},1}\otimes u^{(1)}_{k_{n},2},u^{(2)}_{k_{n},1}
\otimes u^{(2)}_{k_{n},2})$ of $(u^{(1)}_{k_{n}},u^{(2)}_{k_{n}})$ 
up to conjugacy, where
\[
u^{(1)}_{k_{n},1}\equiv S(p_{k_{n}}^{j_{n}})\ ,\ \ 
u^{(2)}_{k_{n},1}\equiv \Omega (p_{k_{n}}^{j_{n}},
\lambda_{k_{n}}^{-1})
\]
for some $1\leq j_{n}\leq i_{k_{n}}$. 
Then it is easy to see that 
\[
(\, {\rm Sp}\left(\,{u^{(1)}_{k_{n},1}}^{p}
{u^{(2)}_{k_{n},2}}^{q}\right)\, 
|\, n\in {\bf N}\, )
\]
 is uniformly distributed in ${\bf T}$. 
This ensures that
\[
(\, {\rm Sp}\left(\,\bigotimes_{k=k_{0}}^{l}{u^{(1)}_{k}}^{p}
{u^{(2)}_{k}}^{q}\right)\, |\, l\geq k_{0}\, )
\]
is also uniformly distributed for all $k_{0}$. 
Hence $\alpha_{(p,q)}=\otimes_{k=1}^{\infty}
Ad\, {u^{(1)}_{k}}^{p}{u^{(2)}_{k}}^{q}$ 
has the Rohlin property by \cite[Lemma\ 5.2]{kishimoto2}. 
\hfill \qed

\vspace{4mm}

To prove Theorem \ref{lab10.2} (2) we introduce an invariant 
as a slight generalization of that in \cite{el2}.

\begin{definition}\label{lab14}
{\rm 
For $n\times n$ unitary matrices $U,\ V$ and $\lambda\in{\bf T}$ 
with $\lambda^{n}=1$, if 
$\| \lambda UV-VU\| <2$ then the closed complex path 
$\gamma (t)={\rm det}((1-t)\lambda UV+tVU)\ \ (t\in [0,1])$ 
does not go through zero. We define $\omega_{\lambda} (U,V)$ 
as the winding number of the path $\gamma$ around zero. 
}
\end{definition}

From $\| \lambda UV-VU\| <2 $ 
we can define 
${\rm log}(\lambda^{-1}VUV^{\ast}U^{\ast})$, 
with  ${\rm log}$  the principal branch of the logarithm. 
As is shown in \cite[Lemma\ 3.1]{exel}
\[
\omega_{\lambda} (U,V)=\frac{1}{2\pi i}{\rm Tr}({\rm log}
(\lambda^{-1}VUV^{\ast}U^{\ast}))
\]
with the nonnormalized trace ${\rm Tr}$
on $M_{n}$.\\

\noindent Proof of Theorem \ref{lab10.2} (2).

Take 
$u^{(1)}_{k},u^{(2)}_{k}\in U(M_{q_{k}})$, 
$\lambda_{k}\in {\bf T}$ as in (\ref{5star}) 
for $\alpha$, 
and $v^{(1)}_{k},v^{(2)}_{k}\in U(M_{q_{k}})$, 
$\mu_{k}\in {\bf T}$ for $\beta$ similarly, i.e.,
\[
\beta_{(p,q)}\lceil M_{q_{k}}=
Ad\, {v^{(1)}_{k}}^{p}{v^{(2)}_{k}}^{q} ,
\]
\[
v^{(1)}_{k}v^{(2)}_{k}=\mu_{k}v^{(2)}_{k}v^{(1)}_{k}\ .
\]
First we show that (2.2) implies (2.1). 
To get a contradiction we assume that outer conjugate 
$\alpha ,\beta$ satisfy $[\alpha ]\not\sim [\beta ]$. 
Outer conjugacy means 
\[
Ad\, W_{i}\circ\alpha_{\xi_{i}}=
\gamma^{-1}\circ\beta_{\xi_{i}}\circ\gamma\ \ 
(i=1,2)
\]
for some unitaries $W_{1},W_{2}\in A$ and 
an automorphism $\gamma$ of $A$. 
For any $\varepsilon >0$ we have a positive integer 
$M$ and unitaries $W_{1}' ,W_{2}'$ in 
$\otimes_{k=1}^{M}M_{q_{k}}$ with 
$\| W_{i}-W_{i}'\| <\varepsilon$  ($i=1,2$) so that
\[
\| Ad\, W_{i}'\circ\alpha_{\xi_{i}}-
\gamma^{-1}\circ\beta_{\xi_{i}}\circ\gamma \| <2\varepsilon\ .
\]
By the assumption there exists a positive integer $K>M$ with 
$\lambda_{K}\neq \mu_{K}$. 
Further take a suffciently large $N>K$ such that
\[
\gamma (M_{q_{K}})\subseteq_{\varepsilon}
M_{q_{1}}\otimes\cdots\otimes M_{q_{N}}\ ,
\]
where $X\subseteq_{\varepsilon}Y$ means that for any 
$x\in X$ there is $y\in Y$ satisfying 
$\| x-y\| \leq\varepsilon\| x\|$. 
Here we use the perturbation theorem 
\cite[Corollary\ 6.8]{christensen}, that is, 
for a sufficiently small $\varepsilon >0$ 
we have a unitary $w_{1}$ in $A$ such that
\[
Ad\, w_{1}\circ\gamma (M_{q_{K}})\subseteq 
M_{q_{1}}\otimes\cdots\otimes M_{q_{N}}
\]
and $\| w_{1}-1\| <28\varepsilon^{\frac{1}{2}}$. 
Set $\gamma_{1}=Ad\, w_{1}\circ\gamma,\ B_{1}
=\gamma_{1}(M_{q_{K}})$ and $B_{2}
=M_{q_{1}}\otimes\cdots
\otimes M_{q_{N}}\cap B_{1}'$. Then 
\[
M_{q_{1}}\otimes\cdots\otimes M_{q_{N}}=
B_{1}\otimes B_{2}\ ,
\]
\[
\gamma_{1}\circ Ad\, W_{i}'\circ\alpha_{\xi_{i}}
\circ\gamma_{1}^{-1}(B_{1})=
B_{1}\ ,
\]
\[ 
\gamma_{1}\circ Ad\, W_{i}'\circ\alpha_{\xi_{i}}
\circ\gamma_{1}^{-1}\lceil B_{1}=Ad\, \gamma_{1} (u^{(i)}_{K})\ 
\]
for $i=1,2$. Since $\| Ad\, W_{1}'
\circ\alpha_{\xi_{1}}-\gamma_{1}^{-1}
\circ \beta_{\xi_{1}}\circ\gamma_{1} \|
\leq C_{1}\varepsilon^{\frac{1}{2}}$ 
for some positive constant $C_{1}$ 
independent of $\varepsilon$ we have 
\[
\beta_{\xi_{1}}(B_{1})
\subseteq_{C_{1}\varepsilon^{\frac{1}{2}}}B_{1}\ .
\]
Noting that $B_{1},\ \beta_{\xi_{1}}(B_{1})
\subseteq M_{q_{1}}\otimes\cdots\otimes M_{q_{N}}$, 
we can use \cite[Corollary\ 6.8]{christensen} again. 
So we have a unitary $w_{2}$ in $M_{q_{1}}
\otimes\cdots\otimes M_{q_{N}}$ such that
\[
Ad\, w_{2}\circ\beta_{\xi_{1}}(B_{1})
\subseteq B_{1}\ ,
\]
\[
\| w_{2}-1\| <C_{2}\varepsilon^{\frac{1}{4}}
\]
for some positive constant $C_{2}$. 
As $M_{q_{1}}\otimes\cdots\otimes M_{q_{N}}$ 
is finite-dimensional and
\[
Ad\, w_{2}\circ\beta_{\xi_{1}}(B_{1})= B_{1}\ ,
\]
\[
Ad\, w_{2}\circ\beta_{\xi_{1}}
(M_{q_{1}}\otimes\cdots\otimes M_{q_{N}})= 
M_{q_{1}}\otimes\cdots\otimes M_{q_{N}}\ ,
\]
we have unitaries $U_{1}$ in $B_{1}$ and $U_{2}$ in $B_{2}$ 
such that 
\[
Ad\, w_{2}\circ\beta_{\xi_{1}}\lceil B_{1}=Ad\, U_{1}\ ,
\]
\[
Ad\, w_{2}\circ\beta_{\xi_{1}}\lceil B_{2}=Ad\, U_{2}. 
\]
For these unitaries we have the following estimates:
\begin{eqnarray*}
\| (Ad\, \gamma_{1}(u^{(1)}_{K})-Ad\, U_{1})\lceil B_{1}\|
&=&
\| (\gamma_{1}\circ Ad\, W_{1}'\circ\alpha_{\xi_{1}}
\circ\gamma_{1}^{-1}-
Ad\, w_{2}\circ\beta_{\xi_{1}})\lceil B_{1}\|\\
&\leq&
\| \gamma_{1}\circ Ad\, W_{1}'\circ\alpha_{\xi_{1}}
\circ\gamma_{1}^{-1}-
Ad\, w_{2}\circ\beta_{\xi_{1}}\|\\
&\leq&
C_{4}\varepsilon^{\frac{1}{4}},
\end{eqnarray*}
\begin{eqnarray*}
\| Ad(U_{1}\otimes U_{2})-\bigotimes_{k=1}^{N}Ad\, v^{(1)}_{k}\|
&\leq&
\| Ad\, w_{2}\circ\beta_{\xi_{1}}-\beta_{\xi_{1}}\|\\
&\leq&
C_{4}\varepsilon^{\frac{1}{4}}
\end{eqnarray*}
for some positive constant $C_{4}$. 
Thus we obtain scalars $\eta_{1},\ \eta_{2}$ 
of ${\bf T}$ such that
\[
\| \gamma_{1}(u^{(1)}_{K})-\eta_{1}U_{1}\| 
<4C_{4}\varepsilon^{\frac{1}{4}}\ ,
\]
\[
\| U_{1}\otimes U_{2}-\eta_{2}
\bigotimes_{k=1}^{N}v^{(1)}_{k}\|
\leq\
4C_{4}\varepsilon^{\frac{1}{4}} .
\]
Consequently we have
\[
\| \gamma_{1}(u^{(1)}_{K})\otimes U_{2}
-\eta\bigotimes_{k=1}^{N}v^{(1)}_{k}\| 
<C_{5}\varepsilon^{\frac{1}{4}}
\]
for some $\eta\in {\bf T}$ and positive constant 
$C_{5}$. Similarly for the direction of $\xi_{2}$ 
we obtain a unitary  $V_{2}$ in $B_{2}$ such that
\[
\| \gamma_{1}(u^{(2)}_{K})
\otimes V_{2}-\zeta\bigotimes_{k=1}^{N}v^{(2)}_{k}\| 
<C_{6}\varepsilon^{\frac{1}{4}}
\]
for some $\zeta\in {\bf T}$ and positive constant 
$C_{6}$. To use an invariant in Definition 
\ref{lab14}, we set $\mu =\prod_{k=1}^{N}\mu_{k}$. 
Then there is a $\lambda\in{\bf T}$ such that 
$\lambda^{q_{1}\cdots q_{N}}=1$ and
\[
|\lambda\mu -1|+2(C_{5}+C_{6})\varepsilon^{\frac{1}{4}}<2\ .
\]
Note that $\omega_{\lambda}$ is invariant under homotopy 
of unitaries for which $\omega_{\lambda}$ is defined. 
From the above estimates  we have
\[
\omega_{\lambda} (\gamma_{1}(u^{(1)}_{K})\otimes U_{2}
,\gamma_{1}(u^{(2)}_{K})\otimes V_{2})=
\omega_{\lambda} (\eta\bigotimes_{k=1}^{N}v^{(1)}_{k},
\zeta\bigotimes_{k=1}^{N}v^{(2)}_{k})
\]
for a sufficiently small $\varepsilon >0$. 
We now evaluate the both sides to get a contradiction. 
Let 
\[
\lambda_{k}={\rm exp}(2\pi i\cdot\frac{s_{k}}{q_{k}}),\ 
\mu_{k}={\rm exp}(2\pi i\cdot\frac{t_{k}}{q_{k}}),
\]
\[
\lambda ={\rm exp}(2\pi i\cdot\frac{s}{q_{1}\cdots q_{N}})
\]
for some $s_{k},t_{k}\in \{0,\ldots ,q_{k}-1\}$ 
and $s\in \{0,\ldots ,(q_{1}\cdots q_{N}-1)\}$. Then
\begin{eqnarray*}
\lefteqn{
\lambda^{-1}(\gamma_{1}(u^{(2)}_{K})\otimes V_{2})
(\gamma_{1}(u^{(1)}_{K})\otimes U_{2})
(\gamma_{1}(u^{(2)}_{K})\otimes V_{2})^{\ast}
(\gamma_{1}(u^{(1)}_{K})\otimes U_{2})^{\ast}}\\
& = & 
1_{q_{K}}\otimes {\rm exp}\left\{ 2\pi i \left(
\frac{-s}{q_{1}\cdots q_{N}}+\frac{-s_{K}}{q_{K}}\right)
\right\}
V_{2}U_{2}V_{2}^{\ast}U_{2}^{\ast}\\
& = &
{\rm exp}\left( 2\pi i\left\{ \bigoplus_{j=1}^{q_{K}}
\left( \left(\frac{-s}{q_{1}\cdots q_{N}}
+\frac{-s_{K}}{q_{K}}\right)
1_{q_{1}\cdots q_{K-1}q_{K+1}\cdots q_{N}}
+H_{2}\right)\right\}\right)
\end{eqnarray*}
for some $H_{2}\in M_{q_{1}\cdots q_{K-1}
q_{K+1}\cdots q_{N}}$ with ${\rm Tr}(H_{2})\in {\bf Z}$.
Thus
\begin{eqnarray*}
\lefteqn{
\omega_{\lambda} (\gamma_{1}(u^{(1)}_{K})\otimes U_{2},
\gamma_{1}(u^{(2)}_{K})\otimes V_{2})}\\
& = &{\rm Tr}\left\{ \bigoplus_{j=1}^{q_{K}}
\left( \left(\frac{-s}{q_{1}\cdots q_{N}}
+\frac{-s_{K}}{q_{K}}\right)1_{q_{1}
\cdots q_{K-1}q_{K+1}\cdots q_{N}}+H_{2}\right)\right\}
\\
& = &
-s-s_{K}q_{1}\cdots q_{K-1}q_{K+1}
\cdots q_{N}+q_{K}{\rm Tr}(H_{2}).
\end{eqnarray*}
On the other hand
\begin{eqnarray*}
\omega_{\lambda} \left(\eta\bigotimes_{k=1}^{N}v^{(1)}_{k},
\zeta\bigotimes_{k=1}^{N}v^{(2)}_{k}\right)
& = &
\left(\frac{-s}{q_{1}\cdots q_{N}}+\sum_{k=1}^{N}
\frac{-t_{k}}{q_{k}}+n\right)\left( \prod_{l=1}^{N}q_{l}\right)\\
& = &
-s-\sum_{k=1}^{N}t_{k}q_{1}\cdots q_{k-1}q_{k+1}
\cdots q_{N}+nq_{1}\cdots q_{N}
\end{eqnarray*}
for some $n\in {\bf Z}$. Therefore
\[
(s_{K}-t_{K})q_{1}\cdots q_{K-1}q_{K+1}\cdots q_{N}
=q_{K}{\rm Tr}(H_{2})+\sum_{\stackrel{\scriptstyle k=1}
{k\neq K}}^{N}\frac{-t_{k}}{q_{k}}+nq_{1}\cdots q_{N}\ .
\]
Noting that $s_{K}-t_{K}\neq 0$ and $s_{K}-t_{K}$ 
is not divided by $q_{K}$, we have a contradiction.

Next we show that (2.1) implies (2.2). We assume 
$[\alpha ]\sim [\beta ]$. Then we may also assume that 
$\lambda_{k}=\mu_{k}$ for any $k$ since inner perturbation 
does not change outer conjugacy classes. If there is a 
$k_{0}\in {\bf N}$ such that $\lambda_{k}=1$ for any 
$k\geq k_{0}$, then we have the result from Theorem \ref{lab10}. 
If there is no such $k_{0}$, we pick up all the $k$'s with 
$\lambda_{k}\neq 1$ and make the subsequence 
$(\, p_{k_{n}}\, |\, n\in {\bf N}\, )$ of 
$(\, p_{k}\, |\, k\in {\bf N}\, )$. Let 
$\lambda_{k}=\exp (2\pi is_{k}/q_{k})$ 
as before. Then for any $N$
\[
\prod_{k=1}^{k_{N}}\lambda_{k}=
\exp\left(2\pi i\cdot\sum_{n=1}^{N}\frac{s_{k_{n}}}{q_{k_{n}}}\right)\ .
\]
By noting that $s_{k_{n}}\neq 0$ and $q_{k_{n}}$'s 
are relatively prime to the each others,  
$\sum_{n=1}^{N}s_{k_{n}}/q_{k_{n}}$ equals to
\[
\frac{S_{N}}{p_{k_{1}}^{j_{1}}\cdots p_{k_{N}}^{j_{N}}}
\]
in the irreducible form for some positive integers 
$j_{1},\ldots ,j_{N}$ and $S_{N}$. Here we apply Lemma \ref{lab11} 
to two pairs 
$(\otimes_{k=1}^{k_{N}}u^{(1)}_{k},\otimes_{k=1}^{k_{N}}u^{(2)}_{k})$, 
$(\otimes_{k=1}^{k_{N}}v^{(1)}_{k},\otimes_{k=1}^{k_{N}}v^{(2)}_{k})$ 
of unitaries. Then for any $\varepsilon >0$, 
if we take a sufficiently large $N$, these pairs are 
almost conjugate i.e., there is a unitary $w_{1}$ in 
$\otimes_{k=1}^{k_{N}}M_{q_{k}}$ such that
\[
\| Ad\, w_{1}\left(\bigotimes_{k=1}^{k_{N}}u^{(i)}_{k}\right)-
\bigotimes_{k=1}^{k_{N}}v^{(i)}_{k}\| <2^{-1}\varepsilon
\]
for $i=1,2$. We adopt the same method for 
$\otimes_{k=k_{N}+1}^{\infty}M_{q_{k}}$ 
and $2^{-1}\varepsilon$ in place of $\varepsilon$. 
Repeating this procedure as in the proof of Theorem 
\ref{lab10}, we have the result. \hfill \qed

\begin{remark}\label{lab13}
{\rm 
Let $A= M_{3}\otimes M_{2^{\infty}}$ 
and let $\omega_{n}={\rm exp}( 2\pi i/n)$ 
for each $n\in {\bf N}$. We define ${\bf Z}^{2}$-actions  
$\alpha ,\ \beta$ on $A$ by
\[
\alpha_{\xi_{1}}=Ad\, \Omega (3,\omega_{3})\otimes
\left(\bigotimes_{k=1}^{\infty}
Ad\, \Omega (2^{k},\omega_{2^{k}})\right)\ ,\ \ 
\alpha_{\xi_{2}}=Ad\, S(3)\otimes\left(\bigotimes_{k=1}^{\infty}
Ad\, S(2^{k})\right)\ ,
\]
\[
\beta_{\xi_{1}}=id_{M_{3}}\otimes
\left(\bigotimes_{k=1}^{\infty}
Ad\, \Omega (2^{k},\omega_{2^{k}})\right)\ ,\ \ 
\beta_{\xi_{1}}=id_{M_{3}}\otimes
\left(\bigotimes_{k=1}^{\infty}Ad\, S(2^{k})\right).
\]
Then the same arguments as in the first part of the above proof 
show that $\alpha,\ \beta$ 
have the Rohlin property and they 
are not approximately conjugarte. However they are clearly 
outer conjugate.
}
\end{remark}

\vspace{4mm}

Let $\{\, q_{k}\, |\, k\in  K\,\}$ be a finite 
or infinite set of prime numbers. 
We next consider product type ${\bf Z}^{2}$-actions 
on the UHF algebra
\[
\bigotimes_{k\in K}M_{q_{k}^{\infty}}\ ,
\]
where $M_{q_{k}^{\infty}}$ is understood 
as $\otimes_{n=1}^{\infty}M_{q_{k}}$\ .

\begin{theorem}\label{lab10.3}
For the above UHF algebra, any two product type 
${\bf Z}^{2}$-actions with the Rohlin property 
are approximately conjugate.
\end{theorem}

\noindent {\it Proof.}\ \ From Theorem \ref{lab10} 
it is enough to prove that for any product type 
${\bf Z}^{2}$-action 
$\alpha$ on $A$ with the Rohlin property 
and $\varepsilon >0$, there exist an automorphism $\gamma$ 
of $A$ and a product type action $\beta$ 
with the Rohlin property such that 
$\alpha\stackrel{\gamma ,\varepsilon}{\approx}\beta$ and $\beta$ 
has the same form as in 
Theorem \ref{lab10}\ \ i.e., there exist a sequence 
$(\, n_{l}\, |\, l\in {\bf N}\, )$ of positive integers 
and sequences 
$(\, v^{(1)}_{l}\, |\, l\in {\bf N}\, ),
\ (\, v^{(2)}_{l}\, |\, l\in {\bf N}\, )$ 
of unitary matrices such that
\[
v^{(1)}_{l},\ v^{(2)}_{l}\in M_{n_{l}}\ ,
\]
\[
v^{(1)}_{l}v^{(2)}_{l}=v^{(1)}_{l}v^{(2)}_{l}\ ,
\]
\[
(A,\beta_{\xi_{1}},\beta_{\xi_{2}})\cong
\left(\bigotimes_{l=1}^{\infty}M_{n_{l}},
\bigotimes_{l=1}^{\infty}Ad\, v^{(1)}_{l},
\bigotimes_{l=1}^{\infty}Ad\, v^{(2)}_{l}\right)\ .
\]

Let $(\, p_{k}\, |\, k\in{\bf N}\, )$ 
be the prime numbers in the increasing order. 
By definition we find a sequence 
$(\, N_{k}\, |\, k\in{\bf N}\, )$ 
of positive integers in such a way that each $k\in{\bf N}$ 
there are nonnegative integers 
$m^{(1)}_{k},\ldots m^{(N_{k})}_{k}$, unitaries 
$u_{k},v_{k}$ and $\lambda_{k}\in{\bf T}$ satisfying
\[
(A,\alpha_{\xi_{1}},\alpha_{\xi_{2}})\cong
\left(\bigotimes_{k=1}^{\infty}
M_{p_{1}^{m^{(1)}_{k}}\cdots p_{N_{k}}^{m^{(N_{k})}_{k}}},
\bigotimes_{k=1}^{\infty}Ad\, u_{k},
\bigotimes_{k=1}^{\infty}Ad\, v_{k}\right) \ ,
\]
\[
u_{k},\ v_{k}\in U(M_{p_{1}^{m^{(1)}_{k}}
\cdots p_{N_{k}}^{m^{(N_{k})}_{k}}}),
\ \ 
u_{k}v_{k}=\lambda_{k}v_{k}u_{k},
\]
\[
m^{(1)}_{k}+\cdots +m^{(N_{k})}_{k}\neq 0\ .
\]
Set $Q_{k}=p_{1}^{m^{(1)}_{k}}\cdots p_{N_{k}}^{m^{(N_{k})}_{k}}$. 
By Remark \ref{lab30}, $\lambda_{k}^{Q_{k}}=1$ for each $k$. 
If $\lambda_{k}=1$ for any $k$, we are done, 
so we assume that $\lambda_{k}\neq 1$ for some 
$k\in {\bf N}$. Take 
$k_{1}\equiv {\rm min}\{\, k\in {\bf N}\, |\, \lambda_{k}\neq 1\, \}$. 
For $n\geq k_{1}+1$ we define
\[
N(n)=\prod_{k=k_{1}+1}^{n}Q_{k},\ \ \ 
\lambda(n)=\prod_{k=k_{1}+1}^{n}\lambda_{k}
\]
and relatively prime integers $M(n)$ and $K(n)$ by
\[
\lambda (n)={\rm exp}\left( 2\pi i\cdot\frac{K(n)}{M(n)}\right)\ ,
\]
\[
M(n)\leq N(n),\ \ \ K(n)\in \{0,\ldots ,M(n)-1\}.
\]
If we set 
$m^{(i)}(n)=\max\{\, m^{(i)}_{k}\, |\, k_{1}+1\leq k\leq n\,\}$ 
($i=1,\ldots ,N_{k_{1}}$) then for each $i$ the exponent 
of the factor $p_{i}$ in $M(n)$ is less than or equal to 
$m^{(i)}(n)$. 
Hence taking a sufficiently large $n$, we can make the exponent 
of the factor $p_{i}$ in $N(n)M(n)^{-1}$ as large as we like 
if $m^{(i)}(n)\neq 0$. In particular 
$N(n)M(n)^{-1}$ is divided by $Q_{k_{1}}^{2}$ and 
$N(n)$ is much larger 
than $M(n)$. 

We want to show that for any $\delta >0$, 
which is much smaller than $\varepsilon$, 
there exist positive integers $n,\ m_{1},\ m_{2}$ 
and unitary matrices 
$U_{i},V_{i}\ (i=1,2,3)$ 
and $W$ such that
\[
n\geq k_{1}+1,\ \ m_{1},m_{2}\geq \delta^{-1},
\]
\[
m_{1}m_{2}=N(n)M(n)^{-1}Q_{k_{1}}^{-2}, 
\]
\[
U_{1},\ V_{1}\in M_{Q_{k_{1}}},\ \ 
U_{2},\ V_{2}\in M_{M(n)Q_{k_{1}}},
\]
\[
U_{3},\ V_{3}\in M_{N(n)M(n)^{-1}Q_{k_{1}}^{-2}},\ \ 
W\in M_{N(n)},
\]
\[
U_{1}V_{1}=\lambda_{k_{1}}^{-1}V_{1}U_{1},\ \ 
U_{3}V_{3}=V_{3}U_{3},
\]
\begin{equation}
{\rm Sp}(U_{3},V_{3})\ is\ (m_{1},m_{2};\delta )-distributed
\label{distrib}
\end{equation}
and
\[
\left(\ \bigotimes_{k=k_{1}+1}^{n}M_{Q_{k}},\ 
\bigotimes_{k=k_{1}+1}^{n}Ad\, u_{k},\ 
\bigotimes_{k=k_{1}+1}^{n}Ad\, v_{k}\ \right)
\stackrel{Ad\, W,2^{-1}\varepsilon}{\approx}
\]
\[
(\ M_{Q_{k_{1}}}\otimes M_{M(n)Q_{k_{1}}}
\otimes M_{N(n)M(n)^{-1}Q_{k_{1}}^{-2}},\ 
\]
\[
Ad(U_{1}\otimes U_{2}\otimes U_{3}),\ 
Ad(V_{1}\otimes V_{2}\otimes V_{3}\ ))\ . 
\]
Suppose that we have shown this statement, then we can 
construct the required $\gamma$ and $\beta$ 
as follows. Set
\[
n_{1}=(Q_{1}Q_{2}\cdots Q_{k_{1}})\cdot 
Q_{k_{1}}\cdot N(n)M(n)^{-1}Q_{k_{1}}^{-2}\ ,
\]
\[
v^{(1)}_{1}
=\left(\bigotimes_{k=1}^{k_{1}}u_{k}\right)
\otimes U_{1}\otimes U_{3}\ ,
\]
\[
v^{(2)}_{1}
=\left(\bigotimes_{k=1}^{k_{1}}v_{k}\right)
\otimes V_{1}\otimes V_{3}\ ,
\]
\[
W_{1}=W\ .
\]
Then $v^{(1)}_{1},\ v^{(2)}_{1}\in U(M_{n_{1}})$. 
Applying the same method to
\[
\left(\ M_{M(n)Q_{k_{1}}}\otimes
\left( \bigotimes_{k=n+1}^{\infty}M_{Q_{k}}\right) ,\ 
Ad\, U_{2}\otimes\left( \bigotimes_{k=n+1}^{\infty}
Ad\, u_{k}\right),\ 
\right. 
\]
\[
\left.
Ad\, V_{2}\otimes
\left( \bigotimes_{k=n+1}^{\infty}Ad\, v_{k}\right)\ 
\right)
\]
and $2^{-1}\varepsilon$ in place of
\[
\left(\ \bigotimes_{k=1}^{\infty}M_{Q_{k}},\ 
\bigotimes_{k=1}^{\infty}Ad\, u_{k},\ 
\bigotimes_{k=1}^{\infty}Ad\, v_{k}\ \right)
\]
and $\varepsilon$, 
we get $n_{2},v^{(1)}_{2},v^{(2)}_{2},W_{2}$. 
Repeating this procedure 
as in the proof of Theorem \ref{lab10}, 
we obtain an automorphism 
$\gamma$ as the infinite product of 
$Ad\, W_{i}$ ($i=1,2,\ldots$) and an action $\beta$ as 
\[
\beta_{\xi_{i}} =\bigotimes_{l=1}^{\infty}Ad\, v^{(i)}_{l}\ \ 
(i=1,2)\ .
\]
This $\beta$ also has 
the Rohlin property due to (\ref{distrib}). 
Hence we have the required $\gamma$ and $\beta$.

Now we show  the remaining  part of the proof, that is, 
the existence of $n,m_{1},m_{2}\in{\bf N}$,  
unitary matrices $U_{i},V_{i}$ ($i=1,2,3$) and $W$ 
satisfying the prescribed conditions. 
By Lemma \ref{lab11} we can decompose (up to conjugacy) 
\[
\left(\ \bigotimes_{k=k_{1}+1}^{n}M_{Q_{k}},\ 
\bigotimes_{k=k_{1}+1}^{n}u_{k},\ 
\bigotimes_{k=k_{1}+1}^{n}v_{k}\ \right)
\]
into
\[
\left( M_{M(n)}\otimes 
M_{N(n)M(n)^{-1}},U^{(n)}_{1}\otimes U^{(n)}_{2},
V^{(n)}_{1}\otimes V^{(n)}_{2}\right)\ ,
\]
where $U^{(n)}_{1}=S(M(n)),\ 
V^{(n)}_{1}=\Omega (M(n),\lambda (n)^{-1})$ 
and $U^{(n)}_{2},V^{(n)}_{1}$ are some commuting unitary matrices. 
Furthermore we see 
\[
\left(\,{\rm Sp}(U^{(n)}_{2},V^{(n)}_{2})\, |\, n\geq k_{1}+1\, \right)
\]
is uniformly distributed. Actually note that 
${\rm Sp}(U^{(n)}_{2},V^{(n)}_{2})$ 
is unique up to piecewise multiples of 
$(\lambda(n)^{k},\lambda(n)^{l})$ for any $k,l\in{\bf N}$. 
So if ${\rm sup}\{\, M(n)\, |\, n\geq k_{1}+1\, \} =\infty$ 
then it is clearly uniformly distributed. If 
${\rm sup}\{\, M(n)\, |\, n\geq k_{1}+1\, \} <\infty$ then 
there is a positive integer $M$ such that 
$M$ is divided by $M(n)$ for any $n\geq k_{1}+1$. Noting that
\[
\left( \bigotimes_{k=k_{1}+1}^{n}u_{k}\right)
\left( \bigotimes_{k=k_{1}+1}^{n}v_{k}\right)=
{\rm exp}\left( 2\pi i\cdot\frac{K(n)}{M(n)}\right)
\left( \bigotimes_{k=k_{1}+1}^{n}v_{k}\right)
\left( \bigotimes_{k=k_{1}+1}^{n}u_{k}\right)
\]
we have
\[
\left( \bigotimes_{k=k_{1}+1}^{n}u_{k}^{M}\right)
\left( \bigotimes_{k=k_{1}+1}^{n}v_{k}^{M}\right)=
\left( \bigotimes_{k=k_{1}+1}^{n}v_{k}^{M}\right)
\left( \bigotimes_{k=k_{1}+1}^{n}u_{k}^{M}\right)
\]
for any $n\geq k_{1}+1$. This implies 
$u_{k}^{M}v_{k}^{M}=v_{k}^{M}u_{k}^{M}$ for any 
$k\geq k_{1}+1$. Hence by Proposition \ref{lab8}
\[
(\ {\rm Sp} \left( \bigotimes_{k=k_{1}+1}^{n}u_{k}^{M},\ 
\bigotimes_{k=k_{1}+1}^{n}v_{k}^{M}\right)
\ |\ n\geq k_{1}+1\ )
\]
is uniformly distributed. Since
\[
{\rm Sp} \left( \bigotimes_{k=k_{1}+1}^{n}u_{k}^{M},\ 
\bigotimes_{k=k_{1}+1}^{n}v_{k}^{M}\right)=
{\rm Sp}\left( 1\otimes {U^{(n)}_{2}}^{M},1
\otimes {V^{(n)}_{2}}^{M}\right)\ ,
\]
it follows that
\[
(\ {\rm Sp} \left( U^{(n)}_{2},\ V^{(n)}_{2}\right)
\ |\ n\geq k_{1}+1\ )
\]
should be uniformly distributed. 
Using this distribution of the joint spectrum, 
we can make for a sufficiently large $n$
\[
\left(
M_{N(n)M(n)^{-1}},U^{(n)}_{2},V^{(n)}_{2}
\right)
\]
close to (up to conjugacy)
\[
\left(
M_{Q_{k_{1}}}\otimes M_{Q_{k_{1}}}
\otimes M_{N(n)M(n)^{-1}Q_{k_{1}}^{-2}},
U_{3}'\otimes U_{4}'\otimes U^{(n)}_{5},
V_{3}'\otimes V_{4}'\otimes V^{(n)}_{5}
\right)
\]
in norm, where 
\[
U_{3}'=V_{4}'=S(Q_{k_{1}}),\ 
U_{4}'=V_{3}'^{\ast}=\Omega (Q_{k_{1}},\lambda_{k_{1}})
\] 
and ${\rm Sp}(U^{(n)}_{5},V^{(n)}_{5})$ is 
$(m_{1},m_{2};\delta )$-distributed  
for some $m_{1},m_{2}\in{\bf N}$ satisfying
\[
m_{1},m_{2}\geq\delta^{-1},\ \ \ 
m_{1}m_{2}=N(n)M(n)^{-1}Q_{k_{1}}^{-2}\ .
\]
Since $U_{4}'V_{4}'=\lambda_{k_{1}}^{-1}V_{4}'U_{4}'$, 
we obtain desired unitaries in such a way that
\[
U_{1}=U_{4}'\ ,\ \ V_{1}=V_{4}'\ ,
\]
\[
U_{2}=U^{(n)}_{1}\otimes U_{3}'\ ,\ \ 
V_{2}=V^{(n)}_{1}\otimes V_{3}'\ ,
\]
\[
U_{3}=U^{(n)}_{5}\ ,\ \ V_{3}=V^{(n)}_{5}\ .
\]
We complete the proof. \hfill \qed

\vspace{4mm}

Now we sum up the results we have shown so far. 
Let $(\, p_{k}\, |\, k\in {\bf N}\, )$ 
be the prime numbers in the increasing order. 
By Glimm's theorem (\cite{glimm}), 
for any UHF algebra $A$ there exists one 
and only one sequence $(\, i_{k}\, |\, k\in {\bf N}\, )$ 
of nonnegative integers or $\infty$ 
such that $A\cong\otimes_{k=1}^{\infty}M_{p_{k}^{i_{k}}}$, 
where $M_{p_{k}^{\infty}}$ 
is understood as $\otimes_{n=1}^{\infty}M_{p_{k}}$.   
Then our classification of ${\bf Z}^{2}$-actions 
on $A$ is as follows:

\begin{theorem}\label{10.4}
Let $A$ be a UHF algebra with the invariant 
$(\, i_{k}\, |\, k\in {\bf N}\, )$ as above.
\noindent {\rm (1)} If 
$\sharp\{\, k\in {\bf N}\, |\, 1\leq i_{k} <\infty\, \}
=\infty$ then there are infinitely many outer conjugacy classes 
of product type ${\bf Z}^{2}$-actions on $A$ 
with the Rohlin property.\\
{\rm (2)} If 
$\sharp\{\, k\in {\bf N}\, |\, 1\leq i_{k} <\infty\, \}<\infty$ 
and $A$ is infinite-dimensional then 
there is one and only one outer conjugacy class 
of product type ${\bf Z}^{2}$-actions 
on $A$ with the Rohlin property.\\
{\rm (3)} If $A$ is finite-dimensional 
then there is no ${\bf Z}^{2}$-action on $A$ 
with the Rohlin property.
\end{theorem}

\vspace{6mm}

\noindent {\bf Acknowledgements}

\vspace{2mm}

The author  would like to express his gratitude 
to Professor Akitaka Kishimoto for suggesting 
this line of the research and for helpful discussions. 
He also would like to thank the referees for the useful 
comments.

\vspace{4mm}

\noindent Department of Mathematics Hokkaido University, 
Sapporo 060 JAPAN\\
E-mail address: h-nakamu@math.sci.hokudai.ac.jp

\end{document}